\documentclass[journal,10pt]{IEEEtran} 

\usepackage{graphics} 
\usepackage{epsfig} 
\usepackage{mathptmx} 
\usepackage{times} 
\usepackage{amsmath} 
\usepackage{amssymb}  
\usepackage{algorithm,algcompatible}
\usepackage{hyperref}
\usepackage{booktabs}
\usepackage{multirow}
\usepackage{siunitx}
\usepackage{graphicx}
\usepackage{epstopdf}

\DeclareGraphicsExtensions{.png,.eps}

\usepackage{natbib}
\setcitestyle{authoryear,open={(},close={)}}

\usepackage{subcaption}

\graphicspath{{./figures/}{./pictures/}}

\title{\LARGE \bf 
	Parameterized Energy-Optimal Regenerative Braking Strategy for Connected and Autonomous Electrified Vehicles: A Real-Time Dynamic Programming Approach
}

\author{Dohee Kim$^{1}$, Jeong Soo Eo$^{1}$, Kwang-Ki K. Kim$^{2*}$
	\thanks{$^{1}$D. Kim and J. S. Eo are with Electrified Systems Control Research Laboratory, R\&D Division, Hyundai Motor Company, Republic of Korea.
		(Email: {\tt doheekim@hyundai.com}; {\tt fineejs@hyundai.com})}
	\thanks{$^{2}$K.-K. K. Kim is with the Department of Electrical Engineering at Inha University, Republic of Korea. (Email: {\tt kwangki.kim@inha.ac.kr})
	}
	\thanks{$^*$Corresponding author: K.-K. K. Kim
	}
}

\begin{document}

\maketitle
\thispagestyle{empty}
\pagestyle{empty}


\begin{abstract}
This paper presents a vehicle speed planning system called the energy-optimal deceleration planning system (EDPS), which aims to maximize energy recuperation of regenerative braking of connected and autonomous electrified vehicles. A recuperation energy-optimal speed profile is computed based on the impending deceleration requirements for turning or stopping at an intersection. This is computed to maximize the regenerative braking energy while satisfying the physical limits of an electrified powertrain. In automated driving, the powertrain of an electrified vehicle can be directly controlled by the vehicle control unit such that it follows the computed optimal speed profile. To obtain smooth optimal deceleration speed profiles, optimal deceleration commands are determined by a parameterized polynomial-based deceleration model that is obtained by regression analyses with real vehicle driving test data. The parameters are dependent on preview information such as residual time and distance as well as target speed. The key design parameter is deceleration time, which determines the deceleration speed profile to satisfy the residual time and distance constraints as well as the target speed requirement. The bounds of deceleration commands corresponding to the physical limits of the powertrain are deduced from realistic deceleration test driving. The state constraints are dynamically updated by considering the anticipated road load and the deceleration preference. For validation and comparisons of the EDPS with different preview distances, driving simulation tests with a virtual road environment and vehicle-to-infrastructure connectivity are presented. It is shown that the longer preview distance in the EDPS, the more energy-recuperation. In comparison with driver-in-the-loop simulation tests, EDPS-based autonomous driving shows improvements in energy recuperation and reduction in trip time.
\end{abstract}


\section{Introduction}
\label{sec1:intro}
In recent years, battery electric vehicles (BEVs) have become prevalent as they improve air quality in urban areas. In addition to BEVs, hybrid electric vehicles (HEVs), plug-in hybrid electric vehicles (PHEVs), and fuel-cell electric vehicles (FCEVs) are the most popular electrified vehicles powered by dual power sources---an engine and an electric motor---for vehicle traction. For vehicle electrification, various new hardware and software developments have been studied in order to reduce greenhouse gas emissions (ecological driving) and improve fuel economy (economic driving). The readers are referred to existing research monographs~\citep{Boulanger2011,Guzzella2013,Ehsani2018}~and references therein for the technical details and history of development of electrified vehicles. 

Energy-efficient or eco-driving technologies can improve the vehicle energy-efficiency for electrified vehicles during traction and have synergistic effects when incorporated with connected and autonomous vehicle (CAV) technology~\citep{Qi2018,vahidi2018energy,Guanetti2018,watzenig2017comprehensive}. Electrified vehicles can generate energy-efficient speed or acceleration profiles by exploiting the look-ahead information obtained from a high-precision map and connectivity via a vehicle-to-everything (V2X) communication network.  Moreover, they can be employed to obtain an optimal energy management system or eco-friendly adaptive cruise control (Eco-ACC)~\citep{Bae2019}. Generally, the use of energy-efficient speed planning or Eco-ACC that utilizes the preceding forecast tends to require that acceleration or cruise conditions be maintained while minimizing deceleration or braking; in addition, those optimal schemes aim mainly to find optimal input sets on the entire route~\citep{Ozatay2014,Wan2016}. 

However, traffic light information is non-linearly linked with future driving circumstances such as traffic congestion and cut-in/out motions of neighboring vehicles. Therefore, 
traffic events during real driving degrade the energy efficiency of optimal solutions considering the entire driving route. Such events may even make a solution invalid owing to violations of traffic regulations. In~\citep{Hao2019}, the author proposed a vehicle trajectory planning algorithm called the eco-approach and departure (EAD) system, which utilized the information of incoming signals, and it was shown that the EAD system was effective for fuel saving and emission reduction at signalized intersections. With the use of real-time SPaT information at signalized intersections, static optimization problems were formulated to generate a vehicle speed trajectory using the parametric optimization of piece-wise trigonometric-linear functions~\citep{Qi2018} and mixed-integer linear programming~\citep{Lu2016}. 

The use of eco-friendly ADAS (Eco-ADAS) to enhance energy efficiency has also been extensively investigated~\citep{Rauh2011,Koch2014}. It is necessary to develop an advanced driver assistance system (ADAS) in order to provide practical energy savings as a semi-autonomous driving concept that can be regarded as Level 2 or 3 autonomous driving, even though CAV technology is expected to provide a progressive perspective in the future. For example, accelerations are made based on the intention of the driver, and when the Eco-ADAS recognizes deceleration events such as turns or red/yellow traffic lights, which explicitly indicate upcoming speed-reduction requirements, the Eco-ADAS can provide energy-efficient deceleration or braking strategies.

To improve energy recovery benefits from deceleration events occurring at irregular but predictable intervals, this paper proposes an energy-optimal deceleration planning system to maximize regenerative energy when the vehicle approaches an upcoming deceleration event. The preview information that can be obtained from CAV technology is employed to decide predictive deceleration parameters such as traffic light phase transition, target speed, and deceleration planning time. For building an optimal control problem (OCP) to maximize regenerative energy, regenerative power for electrified powertrain with negative sign is designed as a cost function wherein the physical limits of the electrified powertrain are explicitly considered.

The proposed EDPS generates a DP-based speed profile in time domain, while an optimal speed value at each node is computed in backward and a time step is used to determine the distance and slope in the spatial domain. To find practically feasible speed candidates over the prescribed planning horizon, state constraints are dynamically updated by considering the road load and the deceleration preference. For determination of deceleration commands containing vehicular deceleration features and driving circumstance information, a practical deceleration model is designed to generate a smooth deceleration profile for which deceleration time is a design variable determined by an optimal search.

To show that the proposed method maximizes energy recuperation performance for connected and automated electrified vehicles, a virtual driving environment is established. The enhancement of energy recovery potential using preview traffic light information is validated in the virtual urban road environment. Various data pre-processing methods are described to effectively utilize the look-ahead information for an upcoming traffic light, and the energy recovery and dynamic performances are comparatively analyzed by varying preview distances, which correspond to locations that receive the preview information in advance. Moreover, the EDPS employing the predictive driving circumstances, including the upcoming traffic light information, is compared to human drivers recognizing traffic lights with their bare eyes. For comparison tests between the EDPS and human drivers, a driving kit set was installed, which implements realistic driving on the virtual urban road environment with 46 traffic lights. The EDPS is implemented on a virtual road with virtual infrastructure, and 10 human drivers operate the driving kit themselves in the same driving conditions. The test results are compared in terms of accumulated energy recovery and a total driving time. The comparison of results indicate that EDPS has the ability to improve the energy recovery performance and reduce the total trip time.

The main contributions of this paper are summarized as follows:

\begin{itemize}

	\item

	We present an optimal control problem to find an energy-optimal speed profile maximizing the recuperated energy obtained from electrified vehicles. For planning a desired speed profile for upcoming events, predictions of geographic conditions and signal timing are necessary. In the present study, preview information available via V2X technology is used to anticipate driving conditions for a partial route with an upcoming deceleration event, and it is integrated into the optimal control problem.

	\item

	To obtain smooth optimal deceleration speed profiles, a parameterized model is used for deceleration commands. The parameters are dependent on preview information such as residual time and distance and target speed. The key design parameter is deceleration time, which determines the deceleration speed profile to satisfy the residual time and distance constraints, as well as the target speed requirement. This parameterized energy-optimal speed planning strategy is particularly useful for reducing the computation time because it does not require any state or input quantization, which is the main difference from existing dynamic programming approaches to energy-optimal speed planning.

	\item

	Instead of imposing numeric values for maximum and minimum deceleration bounds, we use experimental data of real vehicle driving tests for various deceleration scenarios to determine the upper and lower bounds of the deceleration time determining the features of deceleration speed profiles such as peak deceleration and smoothness. Using this approach to set the bounds of deceleration commands makes the proposed parameterized speed planning strategy more practical compared with other optimal control approaches.

	\item

	To demonstrate the effectiveness of the proposed method, energy recovery performances with various preview distances are compared with the driving results of human drivers in virtual driving tests. The comparisons show that the EDPS decreases the total trip time by $3\%\sim10\%$ and increases energy-regeneration efficiency by $16\%\sim130\%$, depending on the preview distances.

\end{itemize}

The paper is organized as follows: Section~\ref{sec2:prob} presents an optimal control problem that devises a speed planning strategy maximizing the regenerative braking performance for energy recuperation. Section~\ref{sec3:deceleration} presents a practical parameterized deceleration model in which the design parameter is determined by solving an OCP that incorporates preview information and physical deceleration limits. Section~\ref{sec4:test} describes a setup of the virtual urban driving test and presents analyses of the proposed EDPS employing different preview distances and comparisons of the energy recuperation performance between the proposed EDPS and human drivers under braking events. Section~\ref{sec5:conclusion} summarizes the main features and contributions of this study.

\section{Problem Formulation for EDPS}
\label{sec2:prob}
%
\begin{figure}[t]
	\centering
	\includegraphics[width=0.875\columnwidth]{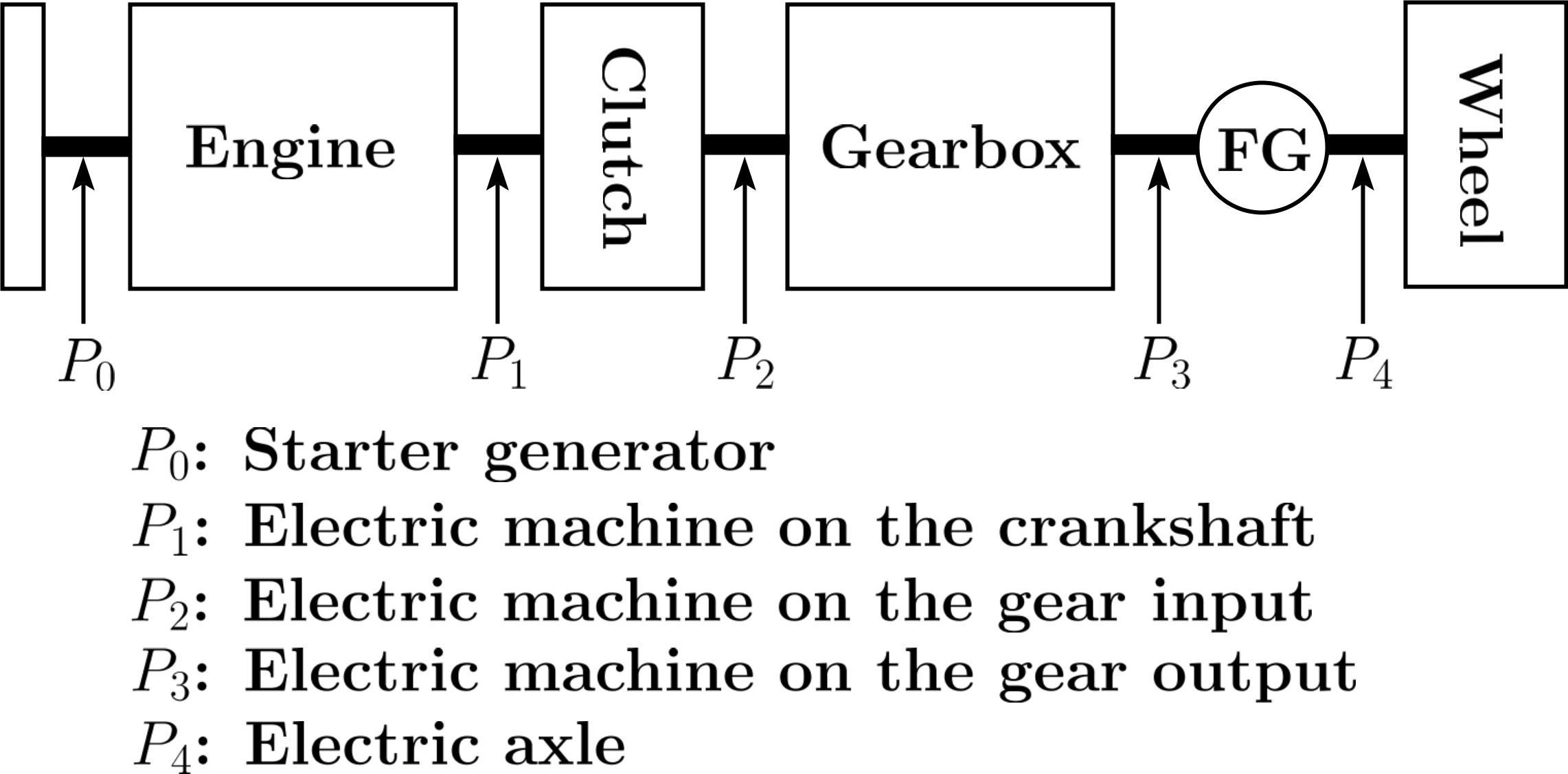}
	\caption{Schematic representation of vehicle electrification architecture: $P_{k}$ denotes the location of an electric machine, and if an electric machine is located on $P_{k}$, it is referred to as $P_{k}$-type electrification. Depending on the locations of electric machines and the combination of power sources, the powertrain model needs to be adjusted in order to reflect changes.}
	\label{fig:elec_arch}
\end{figure}

\begin{figure}[t]
	\centering
	\includegraphics[width=1\columnwidth]{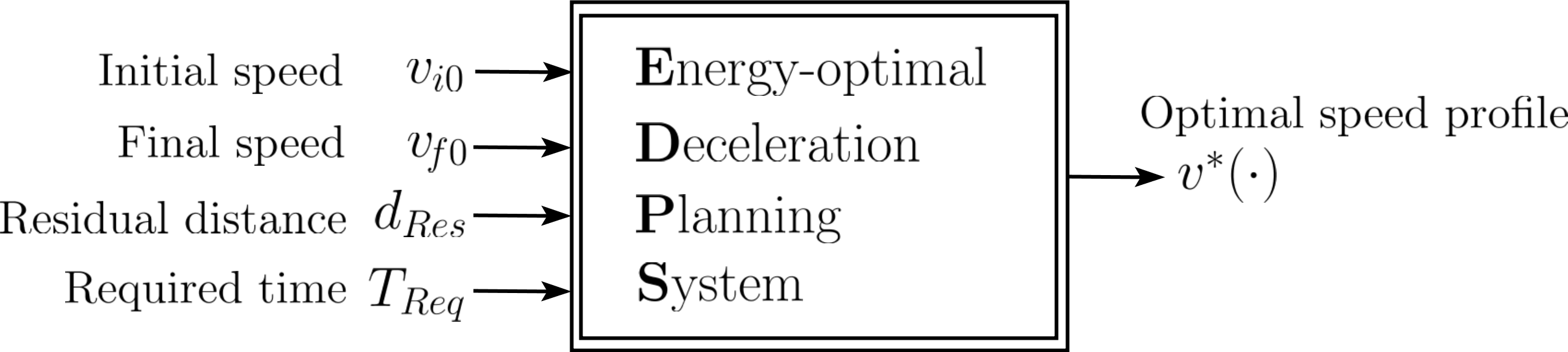}
	\caption{Inputs and outputs of the proposed EDPS for generating a vehicle speed profile that maximizes energy recuperation of regenerative braking.}
	\label{fig:simple_diag}
\end{figure}

For model-based optimal control, this paper considers longitudinal vehicle dynamics of braking and powertrain dynamics of $P_{2}$-type classified
in a vehicle electrification architecture, as shown in Fig.~\ref{fig:elec_arch}. 
To derive the optimal deceleration input generating a speed trajectory to maximize energy regeneration within a deceleration event, EDPS is designed with a DP framework in the time domain. The DP framework assumes that initial and final speeds, remaining distance, and remaining time are available within a predefined distance when approaching a deceleration event. Fig.~\ref{fig:simple_diag} shows the inputs and outputs of the proposed EDPS, where $v_{i0}$ is a current speed of the vehicle, $v_{f0}$ is a required speed at the end of the planning, $d_{Res}$ is a remaining braking-distance, and $T_{Req}$ is a time required for deceleration, which is equal to a planning time of EDPS.

%

\subsection{Optimal Control Problem for Energy Recuperation}
For an OCP formulation, an objective function to be minimized for maximizing energy recuperation in EDPS is given by 
\begin{equation}\label{J}
J=\int_{t_{i}}^{t_{f}}P_{Rgn}\!\left(v\left(t\right),a_{d}\left(t\right),\rho\left(d\left(t\right)\right)\right) dt, 
\end{equation}
where $P_{Rgn}\left(\cdot\right)$ is a systematized regenerative power defined as $F_{Rgn}\left(\cdot\right)v\left(t\right)$, $a_{d}\left(t\right)$ is a deceleration input to be designed, $t_{i}$ and $t_{f}$ are initial and final setting times for the computation of total costs, and $\rho\left(d\left(t\right)\right)$ is a slope information logged for distance $d(t)$ in the spatial domain. It can be explicitly rewritten in terms of the triplet $(v,d,a_d)$ as the product of regenerative braking force induced by an electric motor and vehicle speed:
\begin{equation}\label{eq:sec2:recup_mdl}
\begin{split}
P_{Rgn}\!\left(v, \rho(d), a_{d}\right) &= F_{Rgn} \left(v, \rho(d), a_{d}\right) v. \\
\end{split}
\end{equation}
The deceleration planning problem to maximize energy regeneration can be designed to generate an optimal speed trajectory $v^{*}\left(t\right)$, while finding an acceleration input trajectory $a^{*}_{d}\left(t\right)$ over a {\em fixed} speed-planning time interval $t_{i}\leq t\leq t_{f}$ that minimizes the energy-recuperation performance criterion given in~\eqref{J}.
Note that the smaller the negative value of $J$, the greater is the degree of energy recuperation obtained by regenerative braking.

\subsection{Energy-recuperation Model}
\begin{figure}[t]
	\begin{subfigure}{.25\textwidth}
		\centering
		\includegraphics[width=1.05\textwidth]{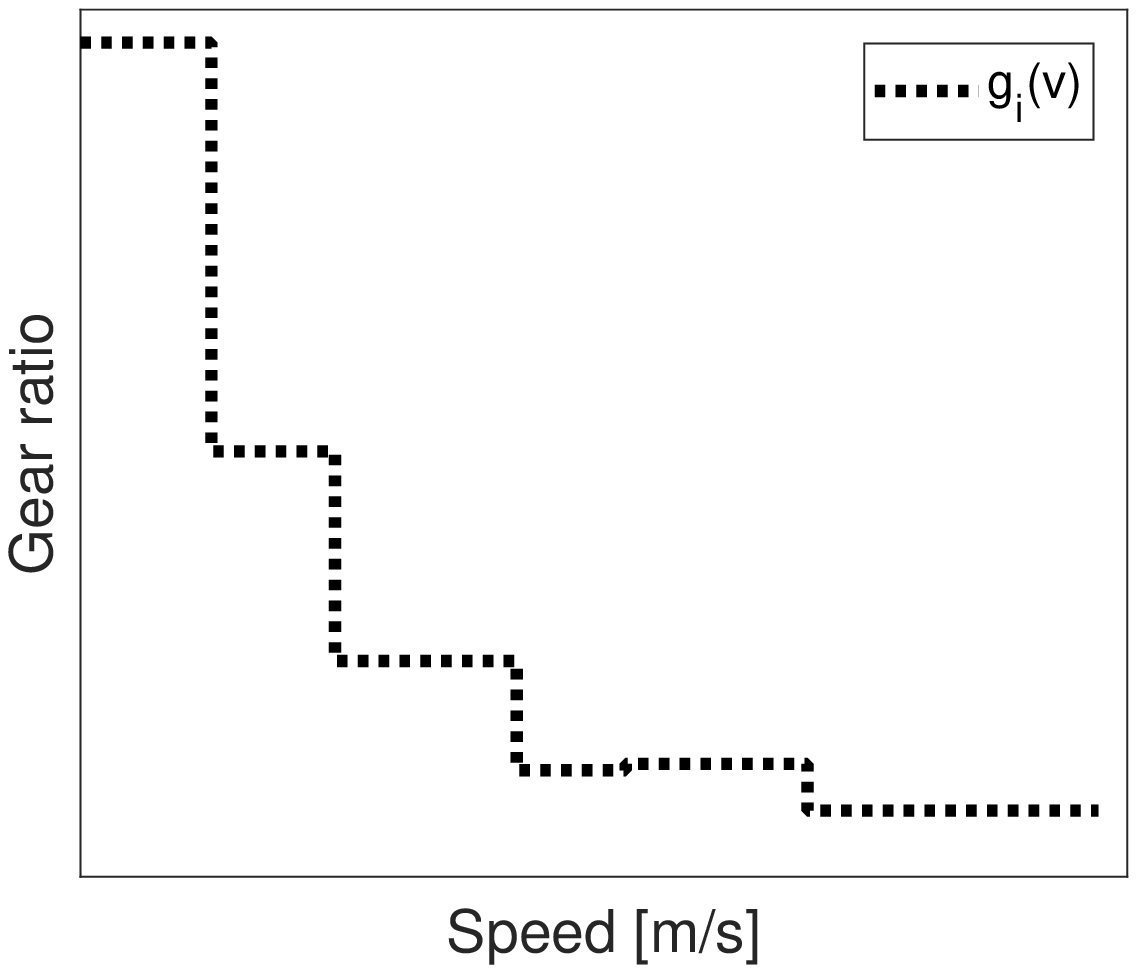}
		\caption{Gear ratio $g_i(v)$}
		\label{fig:Gear_shift}
	\end{subfigure}
	\hspace{-3.75mm}
	\begin{subfigure}{.25\textwidth}
		\centering
		\includegraphics[width=1.05\textwidth]{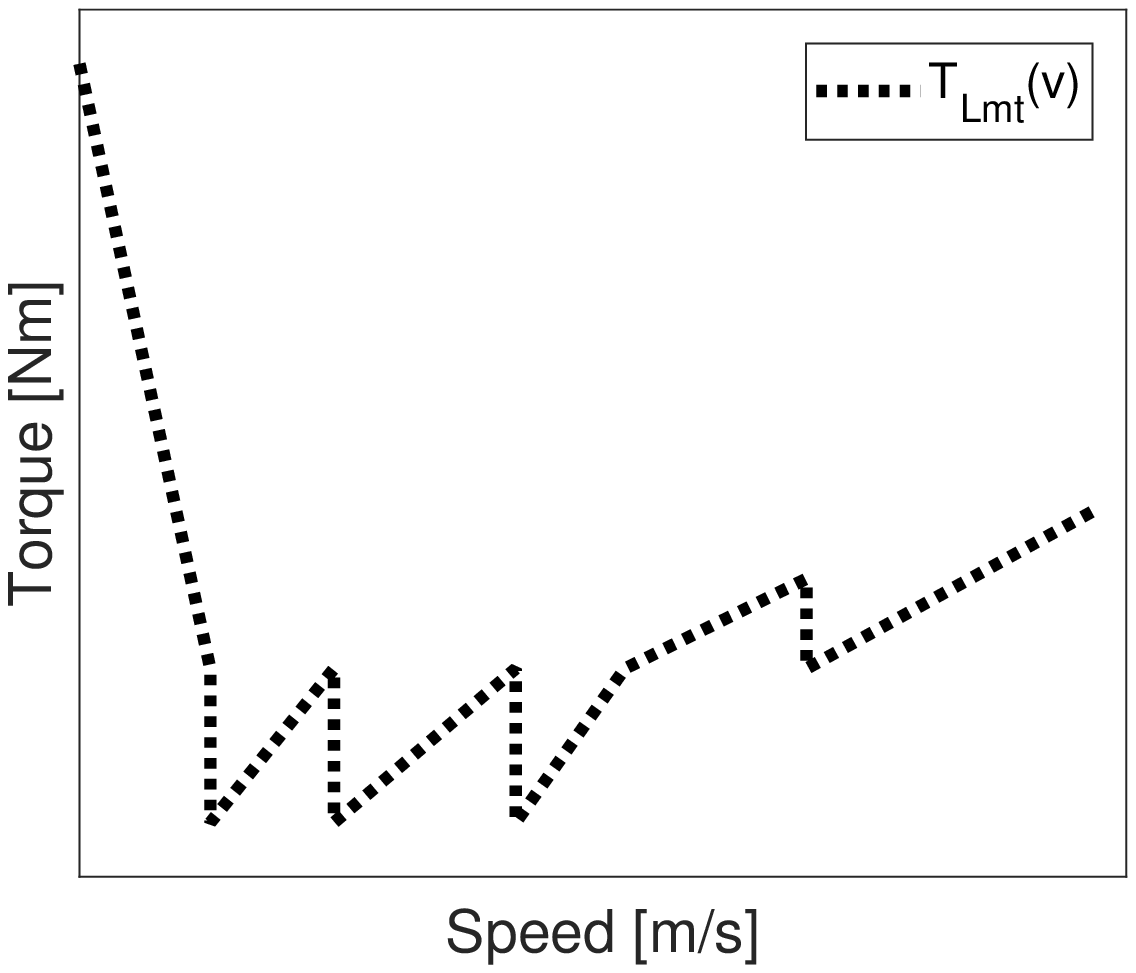}
		\caption{Torque limit $T_{Lmt}(v)$}
		\label{fig:Tq_Lmt}
	\end{subfigure}\vspace{-.5mm}
	\caption{The gear ratios and torque limits of a real-world commercial PHEV using a speed-dependent gear shifting strategy. Owing to proprietary nature of the technical specifications of the car maker, the axes values are not specified.}
	\label{fig:fig}
\end{figure}

\begin{figure}[t!]
	\centering
	\includegraphics[width=1\columnwidth]{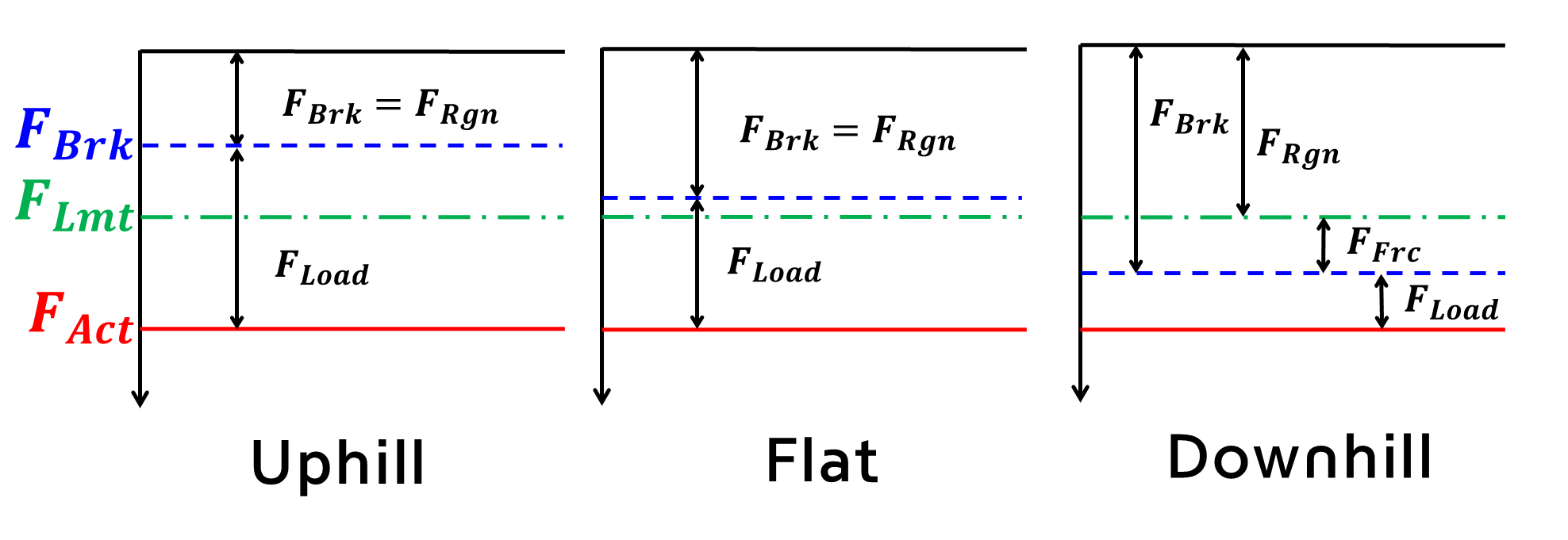}\vspace{-1.5mm}
	\caption{Mechanism snapshots of all of force components in uphill, flat and downhill. In uphill case, $F_{Load}$ exhibits positive increase and is larger than $F_{Load}$ of the flat road, and accordingly, $F_{Brk}$ is reduced and only $F_{Rgn}$ is required for desired deceleration within $F_{Lmt}$. In downhill case, $F_{Load}$ exhibits negative increase and is smaller than $F_{Load}$ of the flat road, and accordingly, $F_{Brk}$ is increased and in addition to $F_{Rgn}\left(=F_{Lmt}\right)$, the mechanical friction force $F_{Frc}$ is required for the desired braking performance.}
	\label{fig:cost_fn}
\end{figure}

The objective of EDPS is to determine the speed profile upon deceleration which minimizes regenerative energy with a negative sign. To accomplish this objective, effective deceleration forces have to be determined in order to contribute to the regenerative deceleration.
By considering the regenerative performance characteristics of the traction motor, the regenerative-braking force can be restricted as
\begin{equation}\label{F_rgn}
F_{Rgn}=\max \left( F_{Brk}, F_{Lmt} \right), 
\end{equation}  
where $F_{Rgn} \leq 0$ is an energy-regenerative deceleration force generated by an electric machine (i.e., a traction motor). Since this paper considers a $P_{2}$-type electrified system, the limited regenerative force of the powertrain, including the traction motor and transmission system, is expressed as
\begin{equation}\label{eq:sec2:motor-force-limit}
F_{Lmt}\left(v\right)=\frac{1}{r_{w}}T_{Lmt}\left(v\right)g_{i}\left(v\right)g_{f},
\end{equation}
where $r_{w}$ is the dynamic wheel radius, $g_{i}\left(v\right)$ indicates the gearbox ratio, which is determined by the gear-shift controller based on the current longitudinal vehicle speed $v$, as depicted in~Fig.~\ref{fig:Gear_shift}, and $g_{f}$ is the final drive ratio. The torque limits of the electric motor in generator-mode are determined by the given motor speed complying with the map given in~Fig.~\ref{fig:Tq_Lmt}, which is obtained by performing quasi-steady state tests that can be mathematically represented as    
\begin{equation}\label{eq:sec2:motor-map}
T_{Lmt}\left(v\right)  = f\left(\omega_{Mot}\left(v\right)\right),
\end{equation}
where $\omega_{Mot}\left(v\right) =\frac{g_i(v) g_f}{ C_{m} r_{w}} v $ is the motor rotation speed in RPM, and
$C_{m} = 2 \pi/ 60$ is the coefficient to convert RPM to ${\rm rad}/{\rm s}$. The quasi-steady state model of the motor-torque limits in rotor speed~\eqref{eq:sec2:motor-map} is used to compute $F_{Lmt}$ by~\eqref{eq:sec2:motor-force-limit}.  
Fig.~\ref{fig:cost_fn} illustrates that $F_{Rgn}$ in the uphill case decreases as the positive slope increases $F_{Load}$, and in the downhill case, $F_{Rgn}$ can be increased by $F_{Load}$ going through the negative slope.
$F_{Rgn} \leq 0$ can be given as
\begin{equation}
F_{Rgn}=F_{Brk}-F_{Frc},
\end{equation}
where $F_{Brk} \leq 0$ is the deceleration force generated from deceleration devices and $F_{Frc} \leq 0$ is a frictional force generated by a hydraulic mechanical braking system. The vehicle braking force $F_{Brk}$ of (\ref{vehicle}) can be presented as
\begin{equation}\label{Force}
F_{Brk}=F_{Load}+F_{Act},
\end{equation} 
where $F_{Load}$ is a road load force and $F_{Act}$ is the inertial force that can be represented as
\begin{equation}\label{vehicle}
F_{Act}=Ma_{d},
\end{equation}
where $M$ is the effective vehicle mass, which is the sum of the curb weight of the vehicle $m$ and the inertia of all of the rotating devices (e.g., motors and engines), and $a_{d}$ is the net acceleration ($a_d > 0$) or deceleration ($a_d \leq 0$) in the longitudinal direction.
Then the longitudinal vehicle dynamics on braking can be written as
\begin{equation}\label{vehicle}
Ma_{d}=F_{Brk}-F_{Load} \,.
\end{equation}
Depending on whether the load force can affect energy-recuperation operation, the road load force $F_{Load}$ can be separated into the following two terms:
\begin{equation}\label{Load}
F_{Load}=F_{Load,\alpha}+F_{Load,\beta}\,.
\end{equation}
The load forces $F_{Load,\alpha}$ and $F_{Load,\beta}$ are defined as 
\begin{equation*}
\begin{split}
F_{Load,\alpha}\left(v,\rho\right)      &=C_{0}\cos\left(\rho\right)+C_{1} v + C_{2} v^{2},\\
F_{Load,\beta}\left(\rho\right) 	      &= Mg\sin\left(\rho\right),
\end{split}
\end{equation*}
where $C_{0}$ and $C_{1}$ are rolling resistance coefficients, and $C_{2}$ is an aerodynamic coefficient. Numeric values of these coefficients are obtained from real vehicle driving tests on normal road surface conditions. The variables $\rho$ and $v$ are the road slope and the longitudinal vehicle speed, respectively. 
To understand the deceleration force that has an effective influence on energy recuperation, $F_{Brk}$ can be decomposed into
\begin{equation}\label{F2}
F_{Brk}=F_{Load,\alpha}+F_{Load,\beta}+F_{Act},
\end{equation}
where $F_{Load,\alpha}$ is a load force that is irrelevant to energy recovery and $F_{Load,\beta}+F_{Act}$ comprises deceleration forces that can contribute to energy recovery. The practical regenerative deceleration force for an electrified vehicle is limited by the generation capability of the electric motor. The feasible regenerative braking force also depends on the powertrain configuration of the electrified vehicle.
%

\subsection{Dynamic Programming Framework}
For speed planning that maximizes energy recuperation, consider the following discrete-time OCP:
\begin{equation}
\begin{split}
\min_{a_{d}(k)\in {\mathcal A}}
&  \  \sum_{k=0}^{N-1} P_{Rgn}\!\left(v(k),a_{d}(k),\rho\left(d(k)\right)\right) \\
\text{ s. t.} \, \, 
& \ v\left(k+1\right) = v(k)+a_{d}(k)\Delta t,\\ \label{DP_J}
& \ d\left(k+1\right) = d(k)+v(k)\Delta t,
\end{split}
\end{equation}
where $\Delta t$ is a sampling-time and $d(k)$ is a distance updated by the determination of $v(k)$.
For an abuse of notation, the notation $v(k)$ is used for $v(k\Delta t)$. The same notation is used for all time-dependent variables.
The constraint set $\mathcal A$ is represented by the following inequalities:
\begin{align}
& d_{\min}(k)	 \leq  d(k) \leq d_{\max} (k) \, , \label{d_r1}\\
& v_{\min}(k)	 \leq  v(k) \leq v_{\max} (k) \, , \label{v_r1}\\
& a_{\min}(k) 	 \leq  a_{d}(k) \leq a_{\max} (k) \, . \label{a_b}  
\end{align}
The lower and upper bounds for speed and deceleration of the set $\mathcal A$ are not fixed but are time-varying or state-dependent.
They are redesigned to integrate varying road load effects and to provide practically feasible operation ranges.

To obtain the solution of the OCP~(\ref{DP_J}),
the Bellman optimality equation for the optimal cost-to-go function $V_{k}$ is given by 
\begin{equation}\label{Bell}
V_{k}\!\left(x(k)\right) = \min_{a_{d}(k)\in {\mathcal A}} \left\{ C(x(k) , a_d(k)) + V_{k+1}\!\left(x(k+1)\right) \right\}, 
\end{equation}
where the state variable is $x(k) = (d(k), v(k))$, and the stage-cost is $C(x(k), a_d(k)) = P_{Rgn} (v(k), a_d(k), \rho(d(k)))$ for $k=N-1, \cdots, 0$.
The terminal condition is given by $V_{N}\! \left( x\left(N\right) \right)=P_{Rgn}\!\left(v(N),0,\rho\left(d\left(N\right)\right)\right)$.

\begin{figure}[t]
	\centering
	\includegraphics[width=0.975\columnwidth]{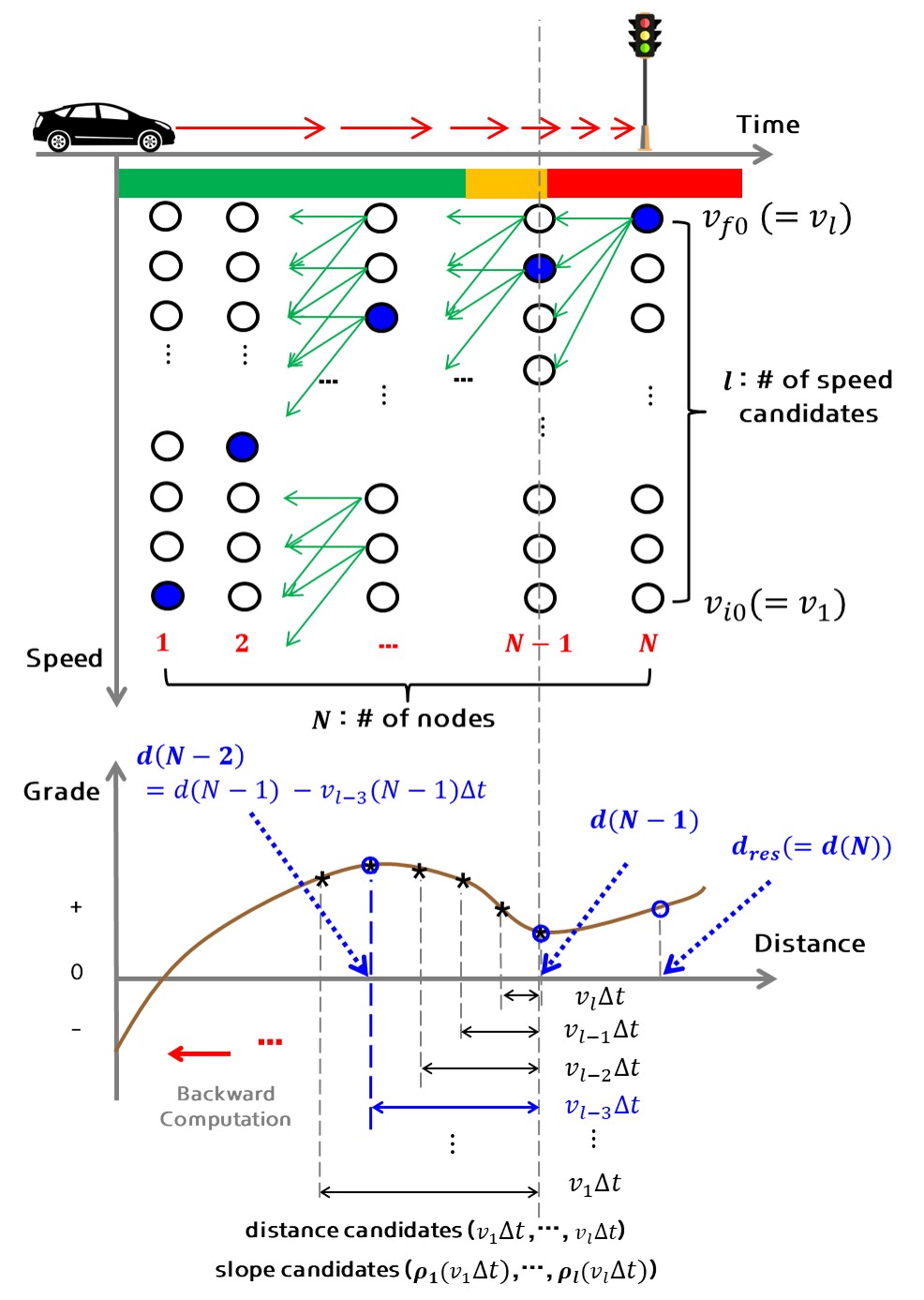}
	\caption{Schematic illustration of the EDPS computation process performed in time and spatial domains: At each time step, the speed profile candidates are determined. From selected speed profile candidates, one can also determine the distance or position profile candidates and then compute the corresponding road slope at time steps over the planning horizon.}
	\label{fig:traffic_DP_desc}
\end{figure}
\begin{figure}[t]
	\centering
	\includegraphics[width=0.999\columnwidth]{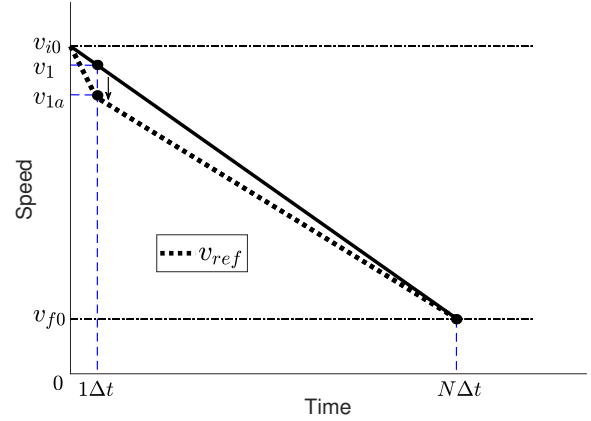}
	\caption{Illustration to build a reference speed profile. The solid line as a monotonically decreasing affine function consisting of two distinct end points $v_{i0}$ and $v_{f0}$ is adjusted to a dashed line by considering the given remaining distance. $v_{1a}$ is a new point determined by the adjustment.}
	\label{fig:spd_ref}
\end{figure}
\begin{figure*}[t]
	\centering
	\includegraphics[width=1.0\linewidth]{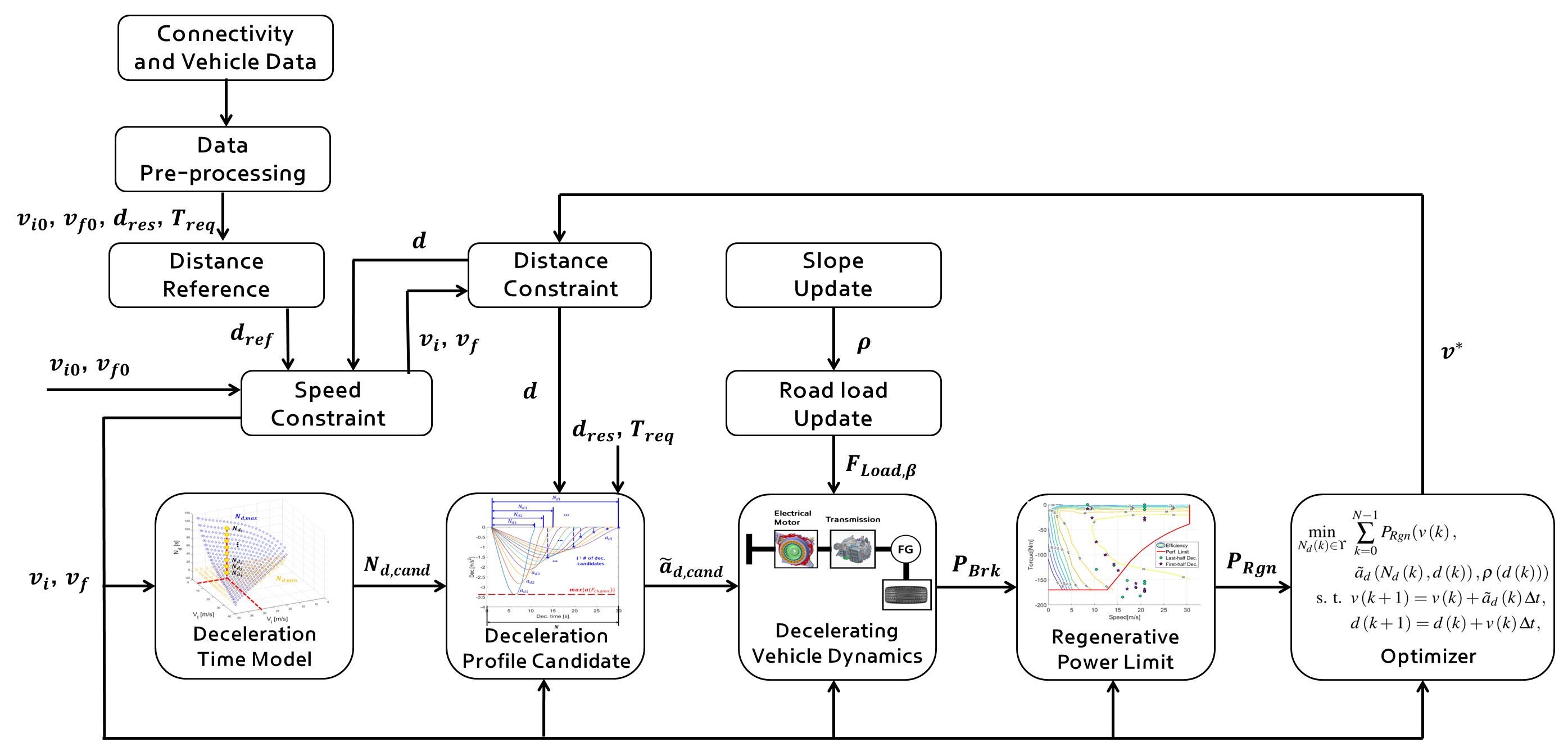}
	\caption{The flow chart of decisions and computations in EDPS.}
	\label{fig:overall_flow}
\end{figure*}
%

\subsection{Computation of Trajectory Candidates}
The system equation for EDPS is computed and updated in the time domain, and the system input to optimize regenerative energy is also determined in the time domain. The representative deceleration events caused by traffic lights are suitable for computation in the time domain because the traffic signal phase and timing (SPaT) restricts the target deceleration time to approach an upcoming traffic light location.
However, because road slope information that dominates road load force is stored in metric form in the spatial domain, optimization in the time domain for utilizing road slope information requires an additional computation method.
A distance in the time domain can be determined using a speed and a specified sampling time step.
As shown in Fig.~\ref{fig:traffic_DP_desc}, an upcoming deceleration event is given within a prescribed remaining distance, and the entire computation node is given as the integer
\begin{equation}
N= \left\lfloor  \frac{T_{Req}}{\Delta t} \right\rfloor  ,
\end{equation}
where $T_{Req}=t_{f} - t_{i}$ is a remaining time for attaining the target deceleration.
When the EDPS generates speed candidates within a speed constraint of (\ref{v_r1}), the speed candidates determine both distance and slope candidates in a series of consecutive time steps, and the specified slope candidates are given by
\begin{equation}\label{slp_cand}
\rho_{1}\left(v_{1}(k)\Delta t\right),\cdots,\rho_{l-1}\left(v_{l-1}(k)\Delta t\right),\rho_{l}\left(v_{l}(k)\Delta t\right),
\end{equation} 
where $l$ is the total number of candidates.
In addition, the optimal speed $v^{*}(k)$ determined at each time step is used to update the remaining distance in a backward computation
\begin{equation}\label{dist_up}
d(k)=d\left(k+1\right)-v^{*}(k)\Delta t,
\end{equation}
where $d(k)$ is the remaining distance updated by the optimal speed at time $k$, and $d\left(N\right)=d_{Res}$ denoted as the remaining distance initially given to EDPS and $k=0,1,\cdots,N-1$.
Fig.~\ref{fig:traffic_DP_desc} shows a schematic for backward computation of the slope candidates and updates of the remaining distances.

\subsection{Dynamic State Constraints}\label{dyn_spd}
To find a set of practically feasible speed candidates at a computation node, the speed constraints in (\ref{v_r1}) need to be dynamically updated by considering the deceleration induced by the road load force, the smooth deceleration for driving preference, and the remaining deceleration distance.
The speed constraints change as the road load forces vary for different driving environments and situations.
In the road load force model of (\ref{Load}), the reference speed profile is designed as the decreasing function of time that complies with the given total deceleration distance requirement.
For dynamic speed constraints, the monotonically decreasing linear function with the two end points $v_{i0}$ and $v_{f0}$ is adjusted in order to satisfy the deceleration distance requirement: 
\begin{equation}
v_{ref}(k)=\frac{1}{N-1}\left(\left(v_{f0}-v_{1a}\right)k+\left(v_{1a}N-v_{f0}\right)\right) \,.
\end{equation}
As illustrated in Fig.~\ref{fig:spd_ref}, the adjusted initial speed considering the remaining distance is given by
\begin{equation}\label{v_1a}
v_{1a}=\frac{2d_{Res}-\left(v_{i0}+v_{f0}\left(N-1\right)\right)\Delta t}{N\Delta t}.    
\end{equation}
The speed change induced by the road load force depending on the reference speed profile is represented by 
\begin{equation}
v_{Load}\left(k+1\right) = v_{Load}(k)- \frac{F_{Load}\left(v_{ref}(k),\rho_{ref}(k)\right)}{M} \Delta t , 
\label{v_load}
\end{equation}
where $\rho_{ref}(k) = \rho( \Delta t \sum_{j=0}^{k} v_{ref}(j) )$ is a slope at time step $k$.
The speed variation determined by (\ref{v_load}) is used to update the upper bound of (\ref{v_r1}) as follows:
\begin{align}\label{ine_vload}
v_{f0} 	&\leq v(k) 	\leq v_{Load}\!(k).
\end{align}

In addition to the road load effects, the speed constraint of \eqref{ine_vload} needs to be dynamically adjusted by comparing the reference distance computed by $v_{ref}$ with the remaining distance determined by $v^{*}$ at each time step.
A distance factor comparing the two distances at time $k$ is defined as
\begin{equation}\label{dist_fac}
\beta(k)=\frac{d_{ref}\left(k+1\right)}{d\left(k+1\right)},
\end{equation}
where $d_{ref}$ is a reference distance that is computed as
\begin{equation}\label{dist_ref}
d_{ref}(k)=d_{ref}\left(k+1\right) - v_{ref}(k)\Delta t \,.
\end{equation}
This distance factor is used to determine the speed range at each time step $k$ as follows:
\begin{align}\label{ine_vdist}
v_{f}\!(k) 	&\leq v(k) 	\leq v_{i}\!(k),
\end{align}
where the upper and lower bounds are respectively given by
\begin{equation} \label{eq:vivf}
\begin{split}
v_{f}(k) &=\max \left(v_{f0},v_{f0}\beta (k)\right) \,, \\
v_{i}(k) &=\min \left(v_{Load}(k),v_{Load}(k)\beta (k)\right) \,,
\end{split}
\end{equation}
which are used to determine the distance range as
\begin{align}\label{ine_d}
d_{f}\!(k) 	&\leq d(k) 	\leq d_{i}\!(k),
\end{align}
where the upper and lower bounds $d_{f}(k)$ and $d_{i}\left(k\right)$ are obtained by $v_{f}(k)$, $v_{i}\left(k\right)$ and the way of distance update in \eqref{dist_ref}.

\section{Practical Design of Deceleration Commands}
\label{sec3:deceleration}

%
\begin{figure}[t]
	\centering
	\includegraphics[width=1.0\columnwidth]{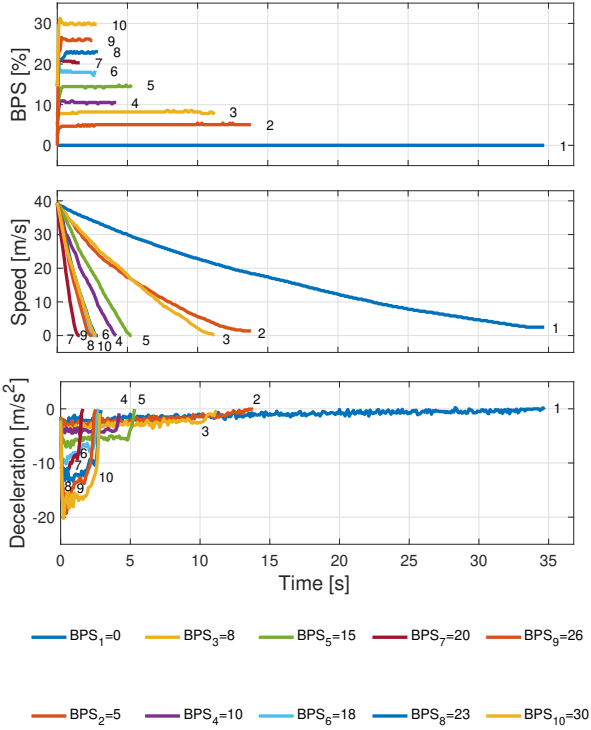}
	\caption{Deceleration test data sets obtained from a real-world commercial PHEV (Hyundai Ioniq) driving test to model deceleration time depending on various combinations of the initial and target speeds. Various BPSs were implemented to evaluate diverse deceleration circumstances.}
	\label{fig:Decel_feature_data}
\end{figure}
\begin{figure}[t]
	\centering
	\includegraphics[width=0.975\columnwidth]{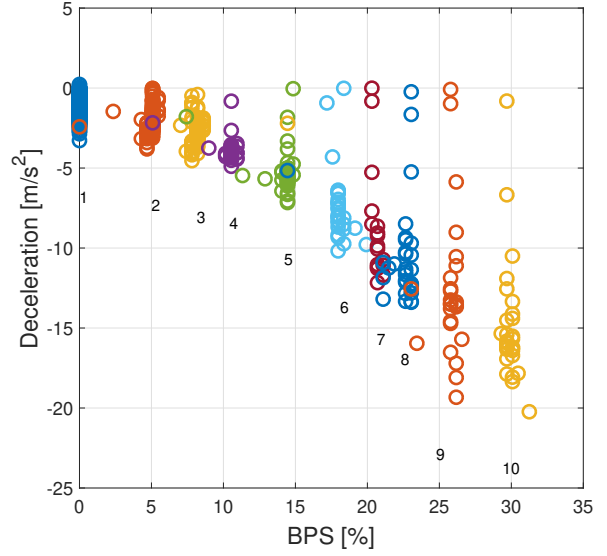}
	\caption{Deceleration over BPS for the entire data set.}
	\label{fig:dec_bps_plot}
\end{figure}
\begin{figure}[t]
	\centering
	\includegraphics[width=0.975\columnwidth]{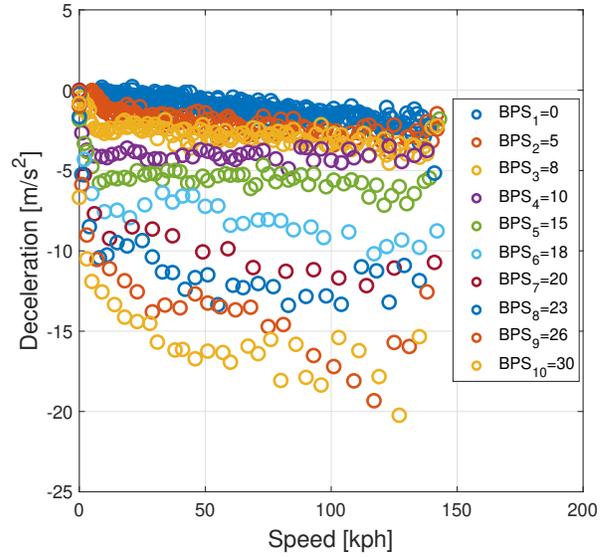}
	\caption{Deceleration over speed for the entire data set.}
	\label{fig:dec_spd_plot}
\end{figure}

This section presents a practical design method for parameterizing and optimizing deceleration commands. As seen in Fig.~\ref{fig:overall_flow}, to determine deceleration commands effectively embracing vehicular deceleration features and driving circumstance information that is given by a connected communication, a practical deceleration model is designed to generate smooth deceleration profiles. This is done by considering deceleration features of a vehicular braking system during the required deceleration time. The deceleration model is presented so that deceleration time can be used as a design variable for determining a deceleration profile. Moreover, main deceleration parameters are determined by analyses of the entire deceleration range of a real vehicle to let the deceleration model include the inherent deceleration features of the vehicle. To reflect the actual deceleration capability of a real-world vehicle, real test data are used to model the actual hardware limits of a market-available vehicle. Fig.~\ref{fig:Decel_feature_data} shows test results of 10 real driving scenarios in which there are 10 different types of braking cases. Upon the modeling and real-driving tests, the deceleration constraints in~\eqref{a_b} are replaced by the deceleration-time constraints that reflect the actual hardware limits of the regenerative and mechanical friction braking systems in a commercial vehicle.

\begin{figure}[t]
	\centering
	\includegraphics[width=0.975\columnwidth]{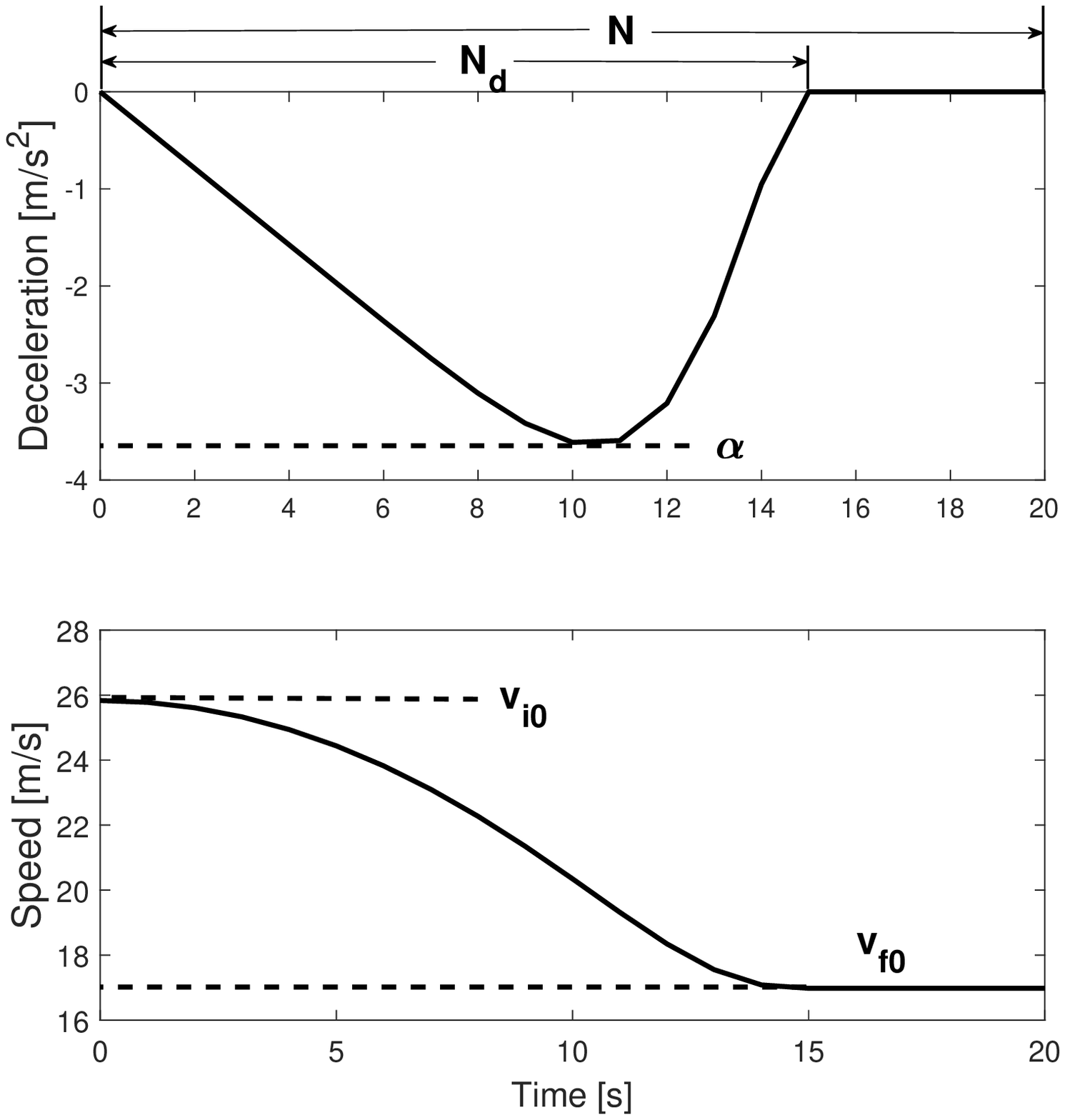}
	\caption{Nominal deceleration and speed profile generated by parameterized deceleration model.}
	\label{fig:base_decel_profile}
\end{figure}
\begin{figure}[t]
	\centering
	\includegraphics[width=0.975\columnwidth]{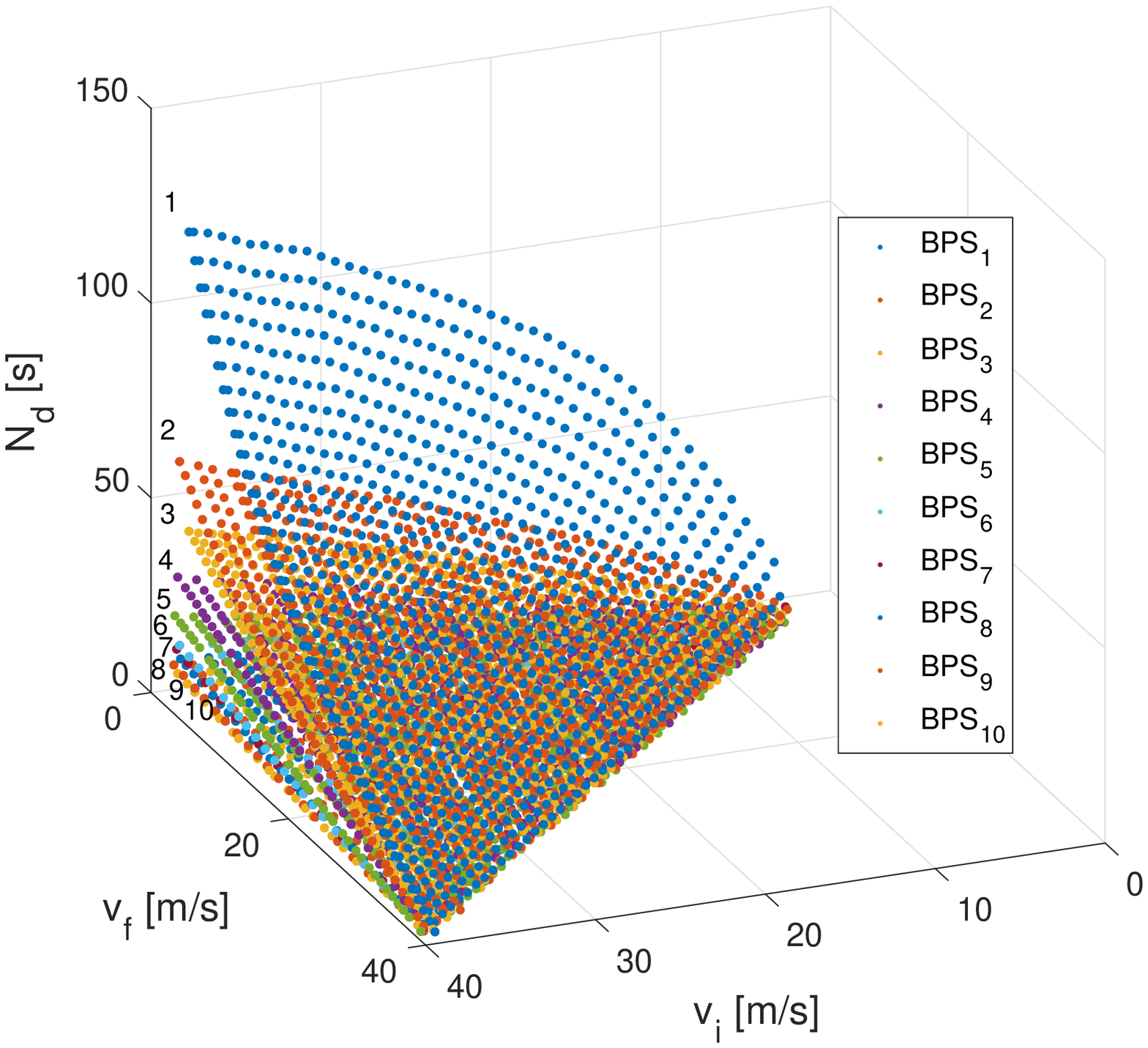}
	\caption{Three-dimensional test data plot of deceleration time over initial speed and final speed for 10 data sets.}
	\label{fig:3d_data}
\end{figure}
%

\subsection{Commercial Vehicle Deceleration Tests}
Practical deceleration limits and features are obtained by analyzing deceleration tests of an actual vehicle. To extract a wide range of deceleration characteristic data, $10$ different braking tests were performed to obtain $10$ different types of vehicle speed profiles ranging from $38.89\: {\rm m/s}$ ($140\: {\rm kph}$) to $0\: {\rm m/s}$ ($0\: {\rm kph}$). In the driving tests, the 10 different sets of braking pedal positions are applied in the commercial vehicle, which are segmented between $0$ (coasting) and $30\:\%$, as shown in Fig.~\ref{fig:Decel_feature_data}. In general, the brake pedal scale (BPS) of $30\:\%$ indicates a rapid deceleration that can be generated by strongly applying the brake pedal by a human driver. It is noticed that the speed and deceleration data are obtained by using a commercial tachometer and the data include the measurement noise. In this paper, a regression model is used to fit these measurement data. 

Deceleration depending on BPS variations in Fig.~\ref{fig:dec_bps_plot} shows a linearly decreasing tendency, except for initial transient parts for each BPS set point. However, a deceleration vs. speed scatter plot given in Fig.~\ref{fig:dec_spd_plot} does not illustrate the linearly decreasing tendency over increasing speed values, but a gradual deceleration over the entire speed domain. For instance, the deceleration can reach about $-16\: {\rm m/s}^{2}$ on $50\: {\rm kph}$ as well as about $-20\: {\rm m/s}^{2}$ on $140\: {\rm kph}$. Such a wide range of deceleration values over each speed segment may result in an unnecessary search of deceleration commands in the optimization process, and may even yield infeasible deceleration at certain time steps. Unless preview information regarding driving circumstances ahead are provided, deceleration commands that can be selected from a wide range of deceleration values can be reasonable control variables to cope with instant variations for uncertain driving circumstances. 
However, provided that preview information such as a target speed and, residual time and distance to the destination are provided, the deceleration time is an effective measure to meaningfully adjust a deceleration profile. As seen in Fig.~\ref{fig:Decel_feature_data}, various deceleration times shown in x-axis can characterize various speed and deceleration profiles.

\subsection{Polynomial-Based Deceleration Model}
To effectively employ preview information provided by connected communication, we propose a smooth deceleration model incorporating the given preview information and the current speed. The proposed deceleration model is inspired from the intelligent driver model proposed in~\citep{Akcelik1987}:
\begin{equation}\label{adip}
\tilde{a}_{d}(N_d(k), d(k))=\left(1-\Delta(k)\right)a_{p}(k)+\Delta(k)a_{pl}(k), 
\end{equation}
where $\Delta(k)$ is a distance ratio presented as
\begin{equation*}
\Delta(k)=w\frac{d(k)}{d_{Res}} \,,
\end{equation*}
where a calibration parameter $w \in [0,1]$ adjusts a deceleration tendency, and $d(k)$ is the remaining distance computed in backward way. In~\eqref{adip}, $a_{p}(k)$ and $a_{pl}(k)$ are defined, respectively, as
\begin{equation}
a_{p}(k)=r\alpha\theta(k)\left[1-\theta^p\!(k)\right]^{2},\label{dip}
\end{equation}
\begin{equation}
a_{pl}(k)=r\alpha\theta_{l}(k)\left[1-\theta^{p}_{l}\!(k)\right]^{2},
\end{equation}
where $p \in {\mathbb R}$, $r = \frac{\left(1+2p\right)^{2+\frac{1}{p}}}{4p^{2}}$, $q = \frac{p^{2}}{\left(2p+2\right)\left(p+2\right)}$,
$\bar a = \frac{v_{f}(k)-v_{i}(k)}{N_{d}(k)}$, and $\alpha = \frac{\bar a}{q}$.
The values of $\alpha$ and $\bar a$ respectively correspond to the maximum and average deceleration with negative signs.
The sequence ratios $\theta(k)$, $\theta_{l}(k)$ are defined by 
\begin{equation}\label{th}
\theta(k)=\frac{k+1}{N_{d}(k)},
\end{equation}
\begin{equation}\label{th_l}
\theta_{l}(k) = \frac{N-k}{N_{d}(k)},
\end{equation}
where $\theta(k),\: \theta_{l}(k)\in \left[ 0,1 \right] \ \mbox{for $k=0,\cdots,N-1$}$, which determine the smoothness and width of the deceleration curve---the larger the deceleration time $N_{d}$, the smoother and wider are the deceleration curves. 
Fig.~\ref{fig:base_decel_profile} shows an example for $a_{p}(k)$, the deceleration profile between $v_{i0}$ and $v_{f0}$.
The deceleration curve can be designed by determining the deceleration time $N_{d}$.
With this parameterized model for shaping deceleration profiles, $N_{d}$ is used as a {\em design} or {\em decision} variable to determine the deceleration command. 
For more details of the polynomial-based deceleration model, please refer to Appendix~\ref{sec:appendix}.
One can observe that $a_{pl}(k)$ flips the shape of $a_{p}(k)$. 
$a_{p}(k)$ with the forward sequence ratio denoted in~\eqref{th} generates a sufficient deceleration in the early phase as seen in Fig.~\ref{fig:base_decel_profile}, whereas $a_{pl}(k)$ with the backward sequence ratio of ~\eqref{th_l} generates a sufficient deceleration in the late phase. As scaling $a_{p}(k)$ and $a_{pl}(k)$ by $\Delta(k)$ that indicates a distance (or position) ratio determined at the current node with a given residual distance, \eqref{adip} can provide a sufficient and smoothing deceleration candidates at each node. From the deceleration data set of Fig.~\ref{fig:Decel_feature_data}, a collection of deceleration times for all combinations of initial and final speeds are illustrated in Fig.~\ref{fig:3d_data}, where the speed gap between initial and final speeds exceeds minimum $2\: {\rm m/s}$. For given initial and final speeds $(v_i(k),v_f(k))$, a candidate set for the deceleration time $N_d(k)$ at each time step $k$ is computed. The $N_d(k)$ map determined by $(v_i(k),v_f(k))$ provides a realistic constraint set including deceleration features of the real vehicle.
\begin{figure}[t]
	\centering
	\includegraphics[width=0.975\columnwidth]{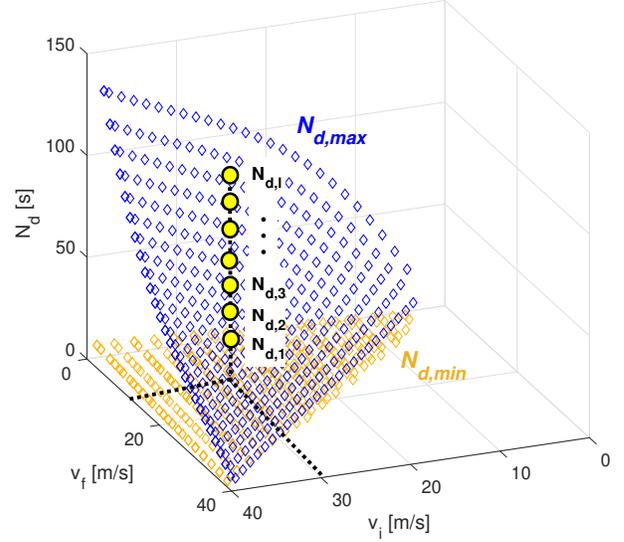}
	\caption{Three-dimensional surface for deceleration time consisting of initial and final speeds (Circles located on the vertical dashed line are extracted by the given $v_{i}$ and $v_{f}$, and are regarded as candidates for $N_{d}$).}
	\label{fig:3D_time_surface}
\end{figure}

Given the speed constraints with the initial and final speed constraints $(v_{i}(k), v_{f}(k))$ and the required planning time $T_{Req}$, the deceleration profile varies with the deceleration time $N_{d}$. As shown in Fig.~\ref{fig:3D_time_surface}, the braking test data set in Fig.~\ref{fig:3d_data} is used to compute the envelope of the feasible deceleration time candidates in terms of $v_{i}$ and $v_{f}$. The two envelopes of Fig.~\ref{fig:3D_time_surface} define the minimum and maximum bounds of $N_{d}(k)$ for a pair of $v_{i}(k)$ and $v_{f}(k)$, respectively. For approximating the boundary surfaces, the following bivariate polynomials are used:
\begin{align}\label{Ndmin}
\begin{split}
N_{d,min} = & \ a_{0}+a_{1}v_{i}+a_{2}v_{i}^{2}+a_{3}v_{i}^{3} + a_{4}v_{f} + a_{5}v_{f}^{2} + a_{6}v_{f}^{3} \, , \\
N_{d,max} = & \ b_{0}+b_{1}v_{i} + b_{2}v_{i}^{2} + b_{3}v_{i}^{3} +b_{4}v_{f} + b_{5}v_{f}^{2} + b_{6}v_{f}^{3} \, ,
\end{split}
\end{align}
where the coefficients are determined by least squares method and the polynomial bases are chosen by trial and error.

\subsection{Practical DP Framework for EDPS}
\begin{figure}[t]
	\centering
	\includegraphics[width=0.975\columnwidth]{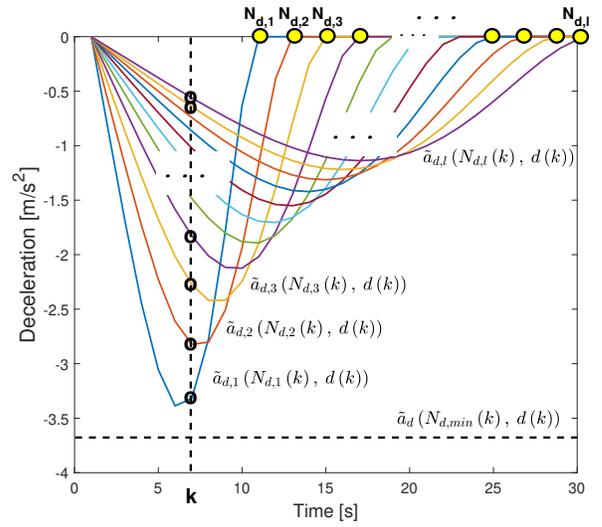}
	\caption{Deceleration profile candidates determined by varying deceleration time candidates and the resulting optimal deceleration inputs denoted by empty circles on each node ($N_{d}$ candidate illustrated as multiple circles generate multiple deceleration profiles as numerous as $N_{d}$ candidates).}
	\label{fig:Decel_cand}
\end{figure}

Deceleration time constrained by the boundary surfaces of~\eqref{Ndmin} can distinguish a deceleration profile from the deceleration model of~\eqref{adip}, and therefore diverse deceleration times generate a variety of deceleration profiles. As depicted in Fig.~\ref{fig:3D_time_surface}, the deceleration time candidates are selected for given speed bounds $\left(v_{i}(k), v_{f}(k)\right)$. Fig.~\ref{fig:Decel_cand} shows the deceleration profiles corresponding to the deceleration time candidates. At every computation node, the multiple candidates, $N_{d,1},\, N_{d,2},\, \cdots,\, N_{d,l}$, generate the multiple deceleration profiles in which the multiple deceleration candidates $\tilde{a}_{d}\!\left(N_{d}(k), d(k) \right)$ are determined, as seen in Fig.~\ref{fig:Decel_cand}. 
To find the optimal input $N_{d}(k)$ to maximize recuperated energy, the DP framework to find the acceleration input of (\ref{DP_J}) is transformed into the framework to find the deceleration time $N_{d}$, the OCP formulation of EDPS is represented by 
\begin{equation}
\begin{split}
\min_{N_{d}(k)\in \Upsilon}
&  \  \sum_{k=0}^{N-1}P_{Rgn}\!\left(v(k),\tilde{a}_{d}\left(N_{d}(k),d(k)\right),\rho\left(d(k)\right)\right) \\
\text{ s. t. } \, \, 
& \ v\left(k+1\right) = v(k)+\tilde{a}_{d}(k)\Delta t,\\ \label{nDP_J}
& \ d\left(k+1\right) = d(k)+v(k)\Delta t,
\end{split}
\end{equation}
where the constraint set $\Upsilon$ is defined by
\begin{align}\label{N_d}
N_{d,min}\!\left(v_{f}(k),v_{i}(k)\right) 	&\leq N_{d}(k) 	\leq N_{d,max}\!\left(v_{f}(k),v_{i}(k)\right), 
\end{align}
where $N_{d,\min}\left(v_{i}(k),v_{f}(k)\right)$ and $N_{d,\max}\left(v_{i}(k),v_{f}(k)\right)$ are determined by (\ref{Ndmin}).
To consider the road slope candidates in order to provide realistic deceleration candidates at each time step, the deceleration time constraint~\eqref{N_d} is modified as
\begin{equation}\label{Nd_adp}
\begin{array}{ll}
N 	\leq N_{d}(k) 	\leq N_{d,max}\!\left(v_{f}(k),v_{i}(k)\right) &\!\!\!\!  \mbox{for } \bar{\rho}(k) \geq 0  \,, \\[1.5mm]
N_{d,min}\!\left(v_{f}(k),v_{i}(k)\right) \leq N_{d}(k) \leq N 	         &\!\!\!\!  \mbox{for }  \bar{\rho}(k) < 0  \,,
\end{array}
\end{equation}
where $\bar\rho(k)$ is the mean of slope candidates~\eqref{slp_cand}~at time step $k$.
The deceleration constraint is given by
\begin{equation}
\tilde{a}_{d}\!\left(N_{d,\min}(k), d(k) \right) \leq \tilde{a}_{d}(N_d(k), d(k))  \leq \tilde{a}_{d}\!\left(N_{d,\max}(k) ,d(k) \right). 
\end{equation}
Using the modified constraints of the deceleration time~\eqref{Nd_adp}, $\tilde{a}_{d}(k)$ can provide mild deceleration candidates with large $N_{d}$ candidates when uphill road is expected to be dominant and provide sharp deceleration candidates with small $N_{d}$ candidates when downhill road is expected to be dominant. 

\begin{figure}[t]
	\centering
	\includegraphics[width=0.975\columnwidth]{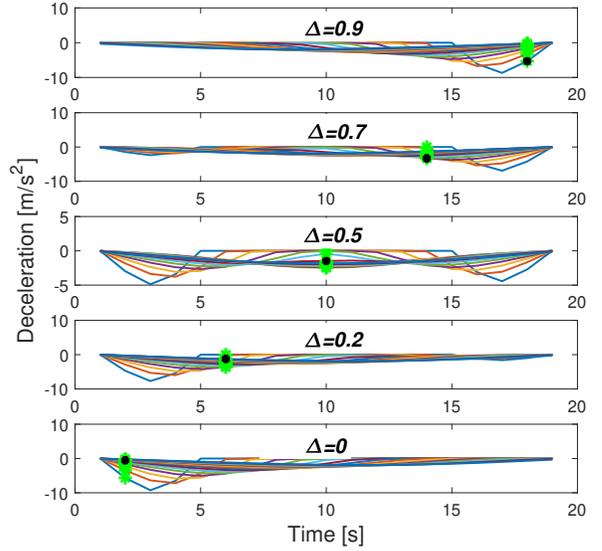}
	\caption{Phased snapshots illustrating the process of selecting an optimal input from multiple deceleration candidates for the deceleration model of~\eqref{adip}.}
	\label{fig:adp_decel}
\end{figure}
Fig.~\ref{fig:adp_decel} shows how the entire set of deceleration candidates changes as the distance ratio $\Delta$ varies from $1$ to $0$.
For $\Delta(k)=1$, $\tilde{a}_{d}(k)=a_{pl}(k)$ and for $\Delta(k)=0$, $\tilde{a}_{d}(k)=a_{p}(k)$.

\begin{figure}[h]
	\centering
	\includegraphics[width=0.975\columnwidth]{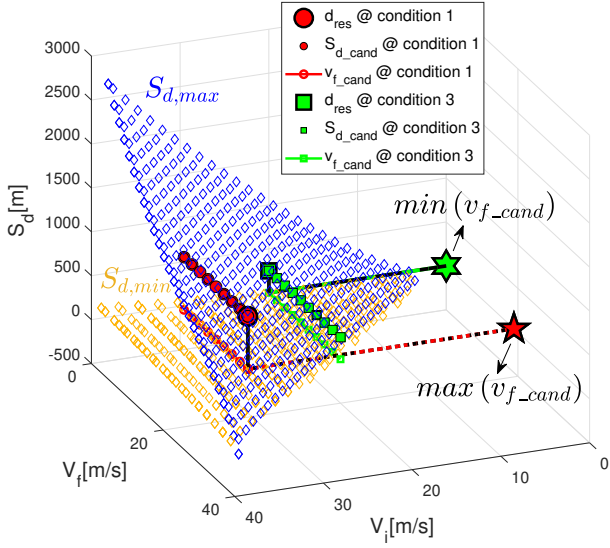}\vspace{-3mm}
	\caption{Two illustrations that describe how to find the target speed. Provided that the remaining distance ($d_{res}$) to an upcoming deceleration event is given, with the deceleration distance model, the distance candidates ($S_{d\_cand}$) varied by target speed candidates ($v_{f\_cand}$) are determined. The transition condition selects the final target speed ($v_{f}$).}
	\label{fig:3d_dist_tgtspd}
\end{figure}
\begin{figure}[t]
	\centering
	\includegraphics[width=0.975\columnwidth]{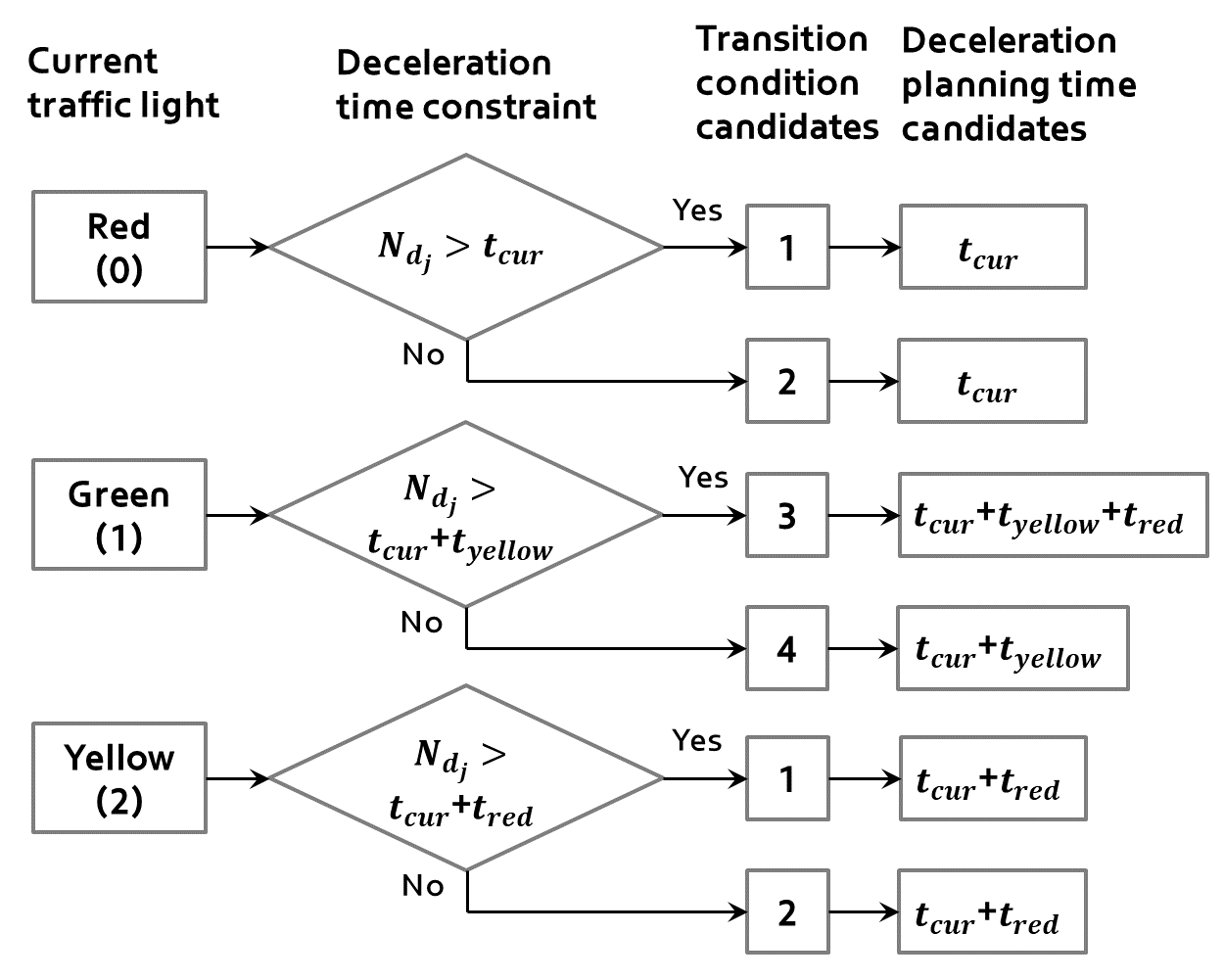}\vspace{-2mm}
	\caption{Flow chart illustrating decision process for deceleration parameters over current traffic light}
	\label{fig:trans_flow_chart}
\end{figure}
%

\subsection{Determination of Deceleration Parameters}
Preview information for upcoming deceleration events are obtained from map-based geographic information and connectivity-based traffic signal information. For map-based deceleration events such as turns and interchange entry/exit, the prescribed speed limit can be a target reference speed. However, deceleration events resulting from traffic lights require an appropriate decision method for determining a specific target speed because traffic light information is not related to the location of the ego-vehicle. Traffic light has its own SPaT which may allow the ego-vehicle to pass by or halt based on the timing at which the deceleration event is detected. Hence, with information of the traffic light ahead, a deceleration condition is determined by the residual distance and a deceleration distance model, which is constructed as
\begin{align}\label{S_d}
S_{d,\min}=\frac{1}{2}  \left( v_{f}+v_{i} \right) {N_{d,\min}} \,, \,
S_{d,\max}=\frac{1}{2}  \left( v_{f}+v_{i} \right) {N_{d,\max}} \,,
\end{align}
where $S_{d,\min}$ and $S_{d,\max}$ are distance boundary surfaces determined based on the initial and final speeds as well as the deceleration time boundary surfaces of (\ref{Ndmin}). The deceleration distance model is an envelope constrained by $v_{i}$, $v_{f}$ and the min/max boundary surfaces, as seen in Fig.~\ref{fig:3d_dist_tgtspd}.

Given the initial speed and residual distance, the modeled deceleration distance candidates are determined to be closest to the given residual distance in the deceleration distance model by the following formulation:
\begin{equation}
\min_{S_{d} \in S_{d\_cand}\left(v_{i}\right)} \left| d_{res} - S_{d} \right| \,.
\end{equation}
The set $S_{d\_cand}\left(v_{i}\right)$ is a partial set of the envelope $\{ S_{d}\left(v_{i},v_{f}\right) : v_{f} \mbox{ is a feasible target speed at final time step}\} $ that is specified as
\begin{align}\label{Sd_cond}
S_{d\_cand}=\left\{S_{d}\left(v_{i},v_{f_1}\right),S_{d}\left(v_{i},v_{f_2}\right),\dots,S_{d}\left(v_{i},v_{f_n}\right)\right\}, 
\end{align}
where $n$ is a total number of multiple distance candidates closest to $d_{res}$, and the target speed set is presented as
\begin{align}\label{vf_cand}
v_{f\_cand}=\left\{v_{f_1},v_{f_2},\dots,v_{f_n}\right\} . 
\end{align}
Since, the envelope of Fig.~\ref{fig:3d_dist_tgtspd} satisfies the linearity condition of (\ref{Ndmin}) and (\ref{S_d}), with (\ref{Sd_cond}), $N_{d\_cand}$ can be linearly transformed into
\begin{align}\label{Nd_cond}
N_{d\_cand}=\left\{N_{d}\left(v_{i},v_{f_1}\right),N_{d}\left(v_{i},v_{f_2}\right),\dots,N_{d}\left(v_{i},v_{f_n}\right)\right\}.
\end{align}

\begin{table}[t]
	\caption{Light phase transition at the start and end of deceleration planning, and the target speed decision for each transition condition.}\label{tab:cond_tab}
	\centering
	\begin{tabular}{|c||@{\,\,}c@{\,\,}|@{\,\,}c@{\,\,}|@{\,\,}c@{\,\,}|}
		\hline
		\textbf{Transition} & \textbf{Light phase} & \textbf{Target} & \textbf{Planning} \\
		\textbf{Condition\:$\left(C_{trans}\right)$} & \textbf{Transition} & \textbf{Speed $\left(v_{f}\right)$} & \textbf{Time $\left(N\right)$}\\
		\hline
		1
		& red~$\rightarrow$~green & 
		$\max{\left(v_{f\_cand}\right)}$ &
		$t_{cur}$ \\
		\hline
		2
		& red~$\rightarrow$~red & $\min{\left(v_{f\_cand}\right)}$ &
		$t_{cur}$ \\
		\hline
		1
		& yellow~$\rightarrow$~green &
		$\max{\left(v_{f\_cand}\right)}$ &
		$t_{cur}+t_{red}$ \\
		\hline
		2
		& yellow~$\rightarrow$~green &
		$\min{\left(v_{f\_cand}\right)}$ &
		$t_{cur}+t_{red}$ \\
		\hline
		3
		& green~$\rightarrow$~red & $\min{\left(v_{f\_cand}\right)}$ &
		$t_{cur}+t_{yellow}+t_{red}$ \\
		\hline
		4
		& green~$\rightarrow$~green & $v_{i}$ &
		$t_{cur}+t_{yellow}$ \\
		\hline		
	\end{tabular} 
\end{table}

To determine the deceleration conditions over the traffic light color for which $0$, $1$, and $2$ denote red, green, and yellow, respectively,
$N_{d_{j}}$, denoted as each component of $N_{d\_cand}$ in~\eqref{Nd_cond}, is compared to the current traffic light duration, yellow light period, and red light period, expressed as $t_{cur}$, $t_{yellow}$, and $t_{red}$, respectively.
In the flow chart, in order to determine the deceleration parameters over a current traffic light shown in Fig.~\ref{fig:trans_flow_chart}, deceleration constraints are investigated by considering the current light color, and then, the transition condition candidates can be denoted as
\begin{equation*}
C_{trans\_cand}=\left\{C_{trans\_cand_1},C_{trans\_cand_2},\dots,C_{trans\_cand_n}\right\} .  
\end{equation*}
The final transition condition for a current traffic light is determined as
\begin{equation}\label{trans_cond}
C_{trans}=\left\lfloor\eta\left(C_{trans\_cand}\right)\right\rceil,
\end{equation}
where $\eta$ represents the median, and $\left\lfloor\cdot\right\rceil$ denotes rounding to the nearest integer. $C_{trans}$ indicates four types of light phase transitions, which predict the traffic light status at the start and end times during deceleration planning, and the light phase transitions are specified in Table~\ref{tab:cond_tab}. As illustrated in Fig.~\ref{fig:3d_dist_tgtspd}, each remaining distance induces each target speed set with the given initial speed and the deceleration distance model. Then, the decision process in Table~\ref{tab:cond_tab} determines the planning time and the target speed denoted as $N$ and $v_{f}$, respectively.

\section{Physics-Based Virtual Urban Driving Tests}\label{sec4:test}
To investigate and verify the energy efficiency of the proposed EDPS, an infrastructure supporting the V2X communication is necessary but is not readily available for experimental tests and fast-track validation. Therefore, the authors have developed a facility of virtual driving tests in which a computer running the EDPS algorithm is connected to a physics-based simulation platform over Ethernet (IEEE 802.3). In the simulation platform, one can customize the test-driving environment parameters such as the road grade and status of traffic lights.

\subsection{Urban Driving Environment Setup}
The virtual driving route is built based on the KOR-NIER Route 1 developed by the National Institute of Environmental Research (NIER) in South Korea. The KOR-NIER Route is a Korean real driving emission (RDE) route, which is in conformance with the trip requirement of the $2^{\rm nd}$ RDE package announced by the European Commission~\citep{Kang_et_al:2017}. The route for the virtual environment was trimmed, and not imitated in the same way, and the traffic SPaT can be arbitrarily set in order to seamlessly implement tests because the environment was established to understand the longitudinal dynamics on the route with the connected infrastructure.

\begin{figure}[t]
	\centering
	\includegraphics[width=0.975\columnwidth]{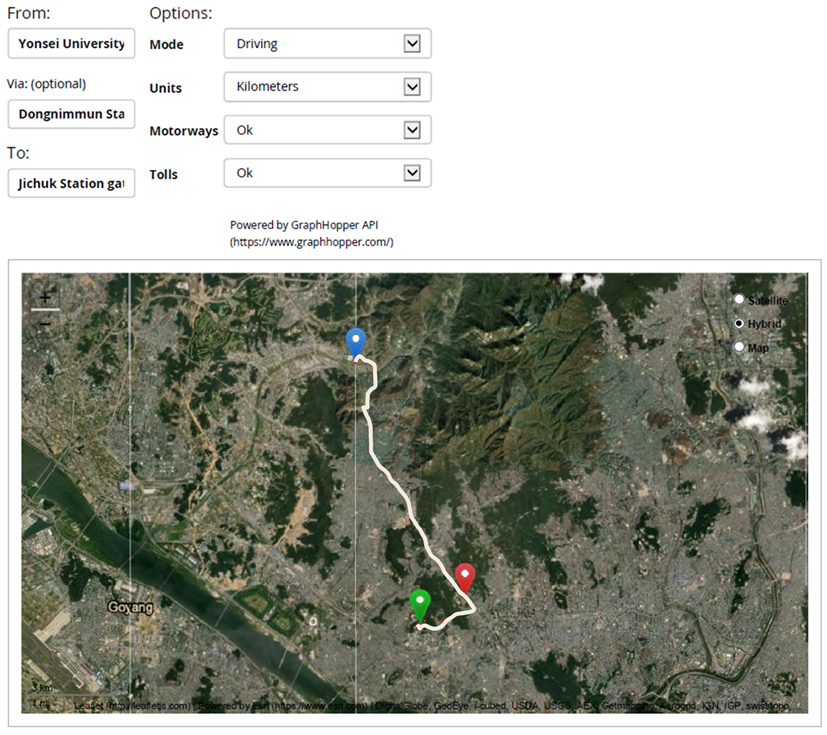}
	\caption{Reference route indicating urban section of KOR-NIER Route 1 for establishing virtual environment.}
	\label{fig:ref_route}
\end{figure}

\begin{figure}[t]
	\centering
	\includegraphics[width=0.975\columnwidth]{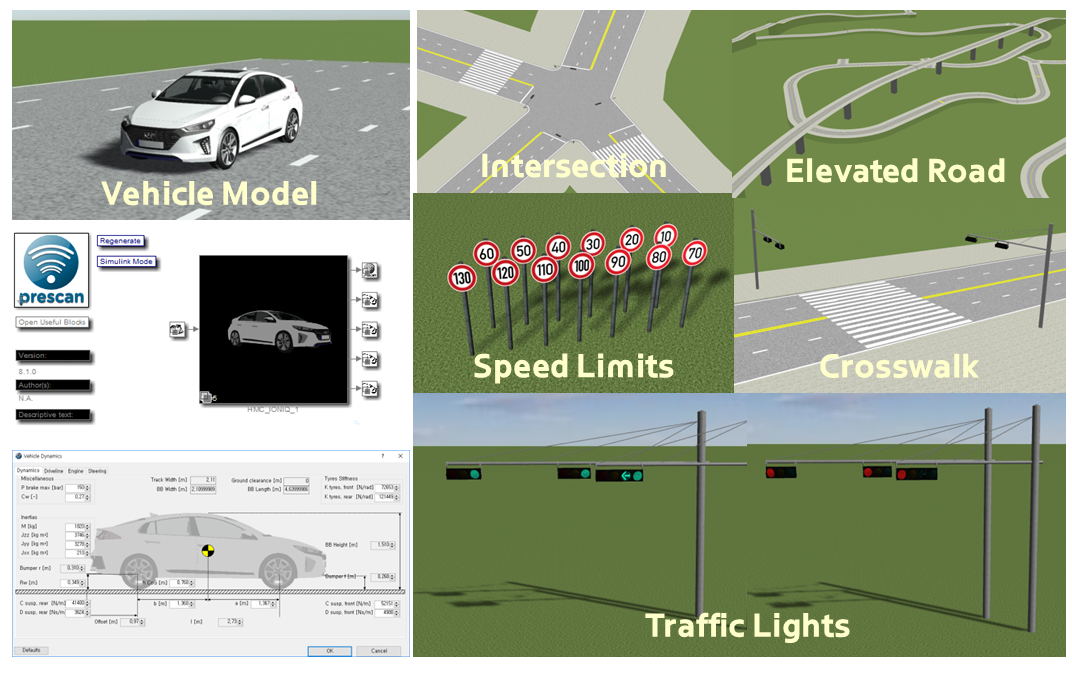}
	\caption{Vehicle model, roads, and hardware equipment that were developed for the virtual environment.}
	\label{fig:virt_equip}
\end{figure}

Fig.~\ref{fig:ref_route} shows an urban section of the KOR-NIER Route 1. The total distance of the route from Yonsei University to Jichuk Station is $12.4\:{\rm km}$,
the road grade varies between $-0.22\:{\rm rad}$ and $0.41\:{\rm rad}$ along the route, and there are $46$ traffic lights. Each traffic light has a period of $25\:{\rm s}$ with $12\:{\rm s}$ for a green light, $10\: {\rm s}$ for a red light, and $3\:{\rm s}$ for a yellow light. In the simulation platform, the vehicle is modeled based on a commercial vehicle, Hyundai Ioniq PHEV. Fig.~\ref{fig:virt_equip} shows a thumbnail sketch for the virtual environment used for the computer-based experiments.

%

\subsection{Connectivity Setting}
\begin{figure}[t]
	\centering
	\includegraphics[width=0.975\columnwidth]{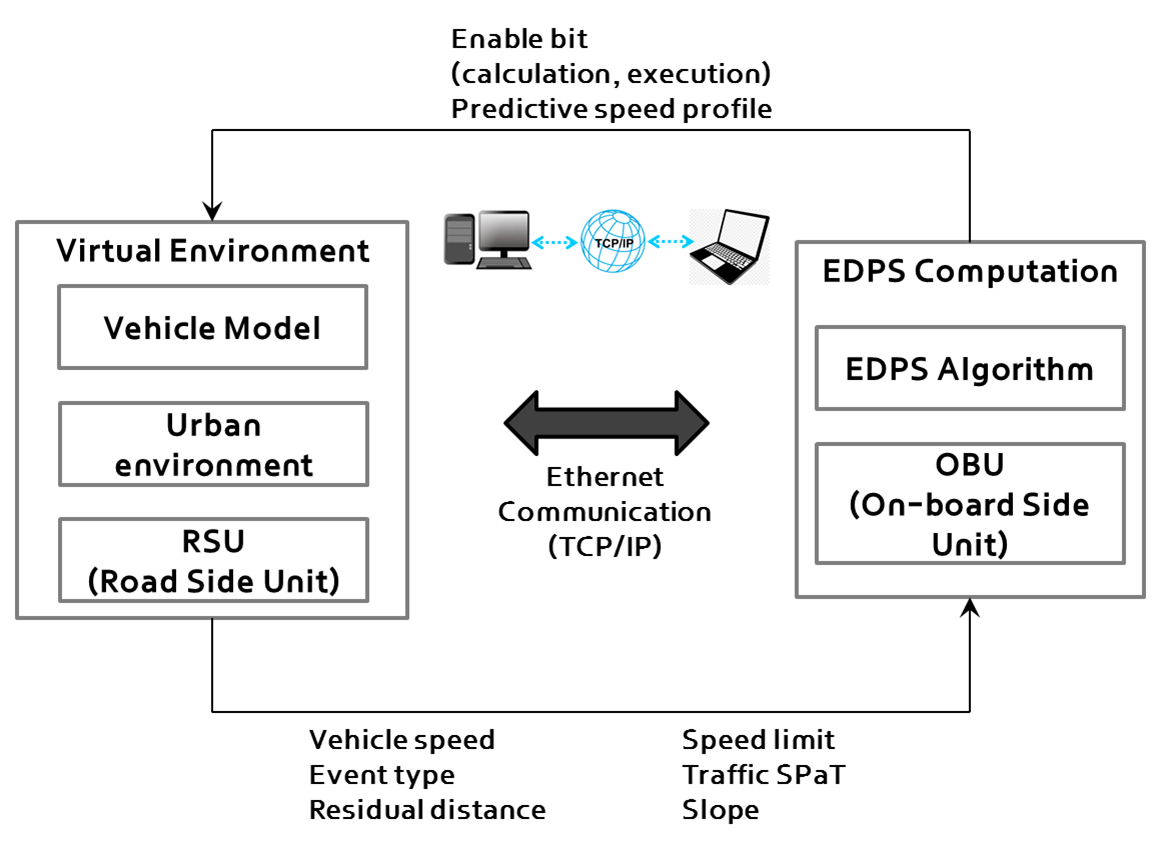}\vspace{-2mm}
	\caption{Communication configuration between computer for simulating virtual driving environment and computation device for computing EDPS.}
	\label{fig:com_config}
\end{figure}

The virtual driving test is performed on two separate computers. As shown in Fig.~\ref{fig:com_config}, the left-side PC (PC 1) includes components depicted in Fig.~\ref{fig:virt_equip} for the virtual environment and generates information, such as the map and traffic SPaT information as well as the vehicle speed. The right-hand side PC (PC 2) computes the EDPS algorithm based on the real-time data given by PC 1. The two PCs are connected by Ethernet using transmission control protocol and internet protocol (TCP/IP).

The experiment scenario is defined as follows:
(i) When the vehicle is within a predefined distance of an upcoming deceleration event, PC 1 transmits the current vehicle speed and virtual environmental information to PC 2;
(ii) PC 2 executes the EDPS with the preview information given from PC 1, and produces a speed profile on deceleration to maximize regenerative energy for the upcoming deceleration event;
(iii) The speed profile is transmitted to the virtual environment of PC 1.
The speed profile can be utilized as a speed reference for the deceleration control of actuators such as the motor or hydraulic braking system. However for the virtual test the profile is executed only on the vehicle assuming  that the actuator precisely tracks the speed reference.

For evaluation of the energy-saving potential that can be achieved by utilizing V2X information such as the traffic SPaT and GPS information, the virtual driving simulation does not consider other vehicles ahead on roads, and it is assumed that the basic driving condition involves cruising at a specified speed. Hence, the vehicle is normally driven in cruise mode with a specified speed until the vehicle meets a deceleration circumstance where the vehicle is within a predefined distance that can receive preview information. When the distance condition for deceleration is met, EDPS plans an energy-optimal speed profile for deceleration and executes the profile for the virtual test. After passing the deceleration event, the vehicle is controlled to linearly accelerate to the specified speed value.

\subsection{EDPS Implementation Results}
%
\begin{table}[t]
	\caption{Summary of EDPS performance over three different preview distance cases.}\label{tab:summary}
	\centering
	\begin{tabular}{|@{}c|c||@{}c@{\,}|@{}c@{\,}|@{}c@{\,}|}
		\hline
		\multicolumn{2}{|c||}{
			$\begin{array}{c} \mbox{Preview} \\ \mbox{distance [meter]} \end{array}$} & $200$  & $150$  & $100$  \\[1mm] \hline
		\multicolumn{2}{|c||}{
			$\begin{array}{c} \mbox{Driving} \\ \mbox{time [sec]} \end{array}$} & $662$ & $639$ & $616$ \\[1mm] \hline
		\multicolumn{2}{|c||}{
			$\begin{array}{c} \mbox{Recuperation} \\ \mbox{energy [Joule]} \end{array}$} & \,$-5.32 \times 10^6$\, &  \,$-5.07 \times 10^6$\,  & \,$-2.74 \times 10^6$\,  \\[1mm] \hline
		\multicolumn{2}{|c||}{
			$\begin{array}{c} \mbox{Average} \\ \mbox{deceleration [m/s$^2$]} \end{array}$} & $-3.79$ & $-3.89$ & $-7.09$ \\[1mm] \hline
		\multicolumn{2}{|c||}{
			$\begin{array}{c} \mbox{Maximum} \\ \mbox{deceleration [m/s$^2$]} \end{array}$} & $-7.47$ & $-7.47$ & $-10.08$ \\[1mm] \hline
		\multirow{4}{*}{$\begin{array}{l@{\!\!}} \mbox{The number of} \\ \mbox{times that it occurs} \\  \mbox{transition} \\ \mbox{conditions} \end{array}$}
		& $C_{trans}=1$ & 5 & 11 & 3 \\ 
		& $C_{trans}=2$ & 6 & 7 & 15 \\
		& $C_{trans}=3$ & 14 & 8 & 3 \\
		& $C_{trans}=4$ & 21 & 20 & 25 \\
		\hline
	\end{tabular} 
\end{table}
%

\subsubsection{EDPS Planning Results}
\begin{figure}[t]
	\centering
	\includegraphics[width=0.975\columnwidth]{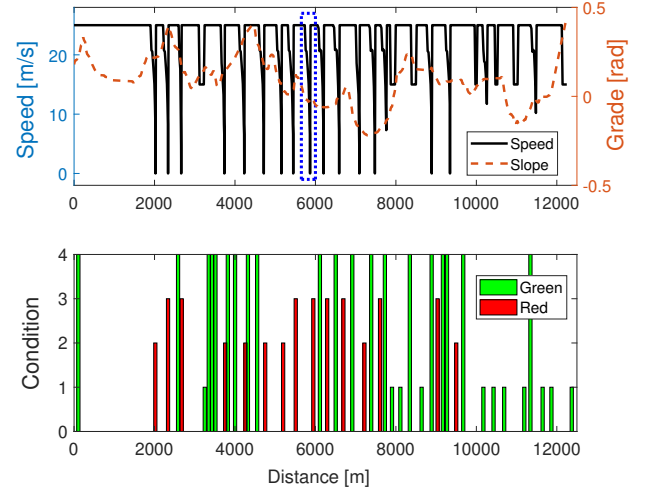}\vspace{-1mm}
	\caption{Regarding a preview distance of $150\:{\rm m}$, speed profile (using left y-axis) and slope profile (using right y-axis) for the entire driving distance at the top of the figure and transition condition sets for each deceleration event at the bottom of the figure.}
	\label{fig:HILS_vel_cond}
\end{figure}

\begin{figure}[h]
	\centering
	\includegraphics[width=0.975\columnwidth]{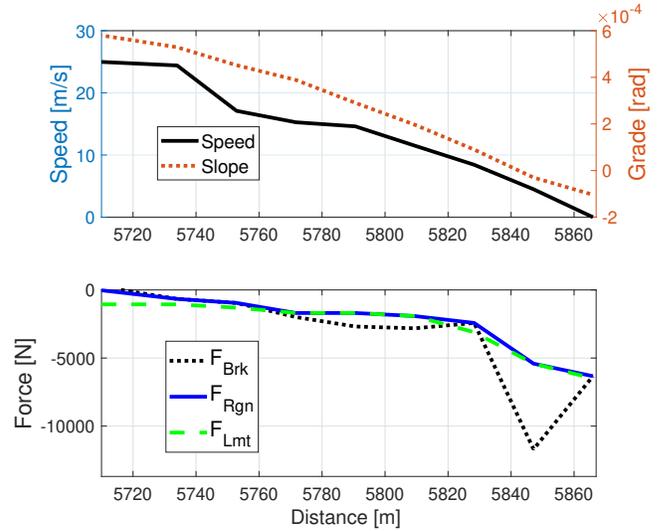}\vspace{-1mm}
	\caption{Speed, slope and force profiles specifying the test results of the rectangle box in Fig.~\ref{fig:HILS_vel_cond}}
	\label{fig:spd_zoom}
\end{figure}

\begin{figure}[t]
	\centering
	\includegraphics[width=0.975\columnwidth]{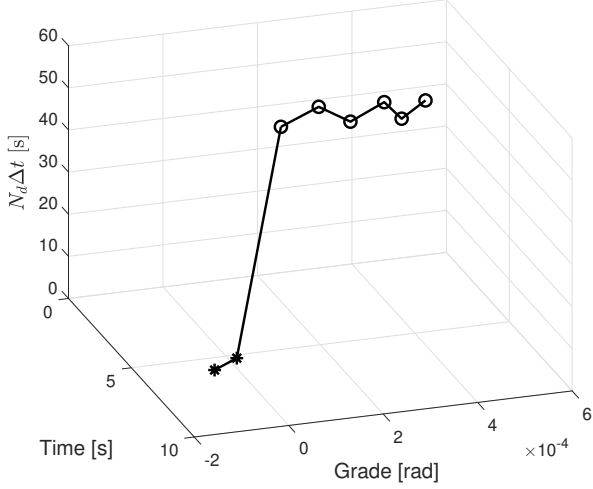}\vspace{-1mm}
	\caption{Optimal $N_{d}$ sequence determined on each node while considering slope information (asterisk: $N_{d}$ on downhill, circle: $N_{d}$ on uphill).}
	\label{fig:Nd_opt}
\end{figure}

\begin{figure}[t]
	\centering
	\includegraphics[width=0.975\columnwidth]{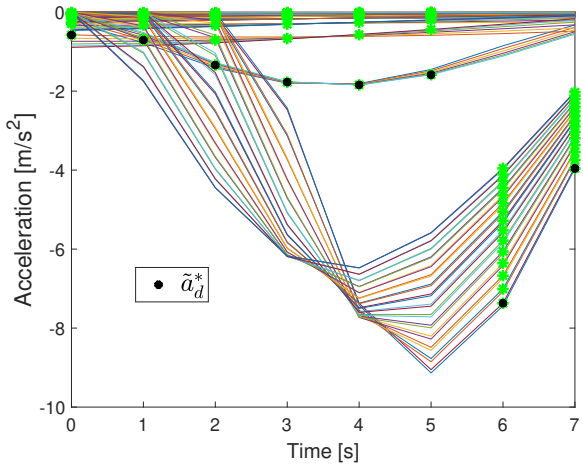}\vspace{-1mm}
	\caption{Overall scenes where a set of deceleration profiles (colorful lines) are generated by a set of deceleration times, a set of deceleration candidates (green asterisks) are determined in each node, and the optimal deceleration (black circle) is selected in each node.}
	\label{fig:acc_hist}
\end{figure}

The top of Fig.~\ref{fig:HILS_vel_cond} displays a vehicle speed profile driven for the entire route, and the bottom of Fig.~\ref{fig:HILS_vel_cond} indicates transition conditions determined by traffic light information transmitted from PC 1 when the vehicle meets the distance condition (specifically, when the vehicle is within $150\: {\rm m}$ from the upcoming traffic lights). In addition, every speed profile on deceleration is computed by EDPS in PC 2. As a result of specifying the rectangle box in the top of Fig.~\ref{fig:HILS_vel_cond} to analyze how the EDPS exploits preview road information, Fig.~\ref{fig:spd_zoom} shows speed and slope profiles with detailed x-axis on the top and the force profiles at the bottom. At the bottom of Fig.~\ref{fig:spd_zoom}, the proposed EDPS generates smaller deceleration forces on the uphill until increasing the deceleration forces on the downhill while $F_{Lmt}$ is determined by considering the motor generation limit and gear box ratio.
Fig.~\ref{fig:Nd_opt} illustrates an optimal deceleration time profile, $N_{d}^{*}(k)$, dependent on the slope variation for the rectangle box case because the deceleration time candidates
denoted as $N_{d}(k)$ reflect the slope information of each node, which are constrained by the inequality given in~\eqref{Nd_adp}.
$N_{d}(k)$ can yield multiple deceleration profile candidates throughout the planning duration using the adapted deceleration model of~\eqref{adip}, as illustrated in Fig.~\ref{fig:acc_hist}, which highlights the overall selection process and results.

Once the distance from a traffic light and its SPaT are provided, the EDPS determines the transition condition in terms of the upcoming traffic SPaT given in~\eqref{trans_cond} and Fig.~\ref{fig:trans_flow_chart} and the target speed is also determined by (\ref{vf_cand}) and Table~\ref{tab:cond_tab}. In the test case with a preview distance of $150\: {\rm m}$ to upcoming traffic lights, the total driving time is $639\: {\rm s}$ for the entire route of $12.4\: {\rm km}$, and the accumulated recuperation energy is $5.07 \times 10^6\: {\rm J}$. To compute the accumulated recuperation energy, $P^{*}_{Rgn}$ consisting of optimal solutions of~\eqref{nDP_J} in the structure of~\eqref{eq:sec2:recup_mdl} is integrated for each deceleration event, and each recuperation energy is accumulated for all deceleration events.
%

\subsubsection{Performance Analysis}
\begin{figure}[t]
	\centering
	\includegraphics[width=0.975\columnwidth]{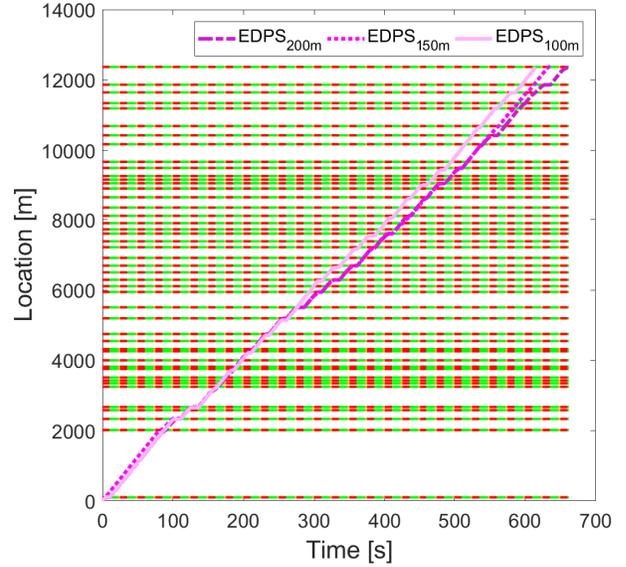}\vspace{-1mm}
	\caption{Location profiles guided by EDPS using $200\:{\rm m}$, $150\:{\rm m}$, and $100\:{\rm m}$ preview distances over driving time and traffic SPaT.}
	\label{fig:EDPS_loc_time}
\end{figure}

To investigate the EDPS performance relying on the distance between an ego-vehicle and upcoming traffic lights, the EDPS performance for various preview distances is tested. Regarding three preview distances of $200\: {\rm m}$, $150\:{\rm m}$ and $100\:{\rm m}$, each EDPS test is executed. Fig.~\ref{fig:EDPS_loc_time} shows the locations of vehicles guided by EDPS using three different preview distances over driving time and the traffic SPaT. The total driving time is measured as $662\:{\rm s}$, $639\:{\rm s}$, and $616\:{\rm s}$, respectively. The longer the preview distance, the longer the deceleration time because at longer preview distance the vehicle decelerates early and steadily. Hence, the EDPS with a preview distance of $200\:{\rm m}$ the most frequently experiences the $3^{\rm rd}$ transition condition denoted as $C_{trans}=3$, which retains the largest energy recovery potential.

A summary of the EDPS performance for each preview distance is provided in Table~\ref{tab:summary}. The EDPS with the $100\:{\rm m}$ preview distance the most frequently belongs to the $4^{\rm th}$ transition condition which passes upcoming traffic lights without any deceleration, and the EDPS shows the least energy recovery but arrives at the final destination the most rapidly. In addition, the $100\:{\rm m}$ test case shows the largest maximum and mean deceleration to halt the vehicle during the short planning time given by the frequent transition conditions of $C_{trans}=2$.

\begin{figure}[t]
	\centering
	\includegraphics[width=0.975\columnwidth]{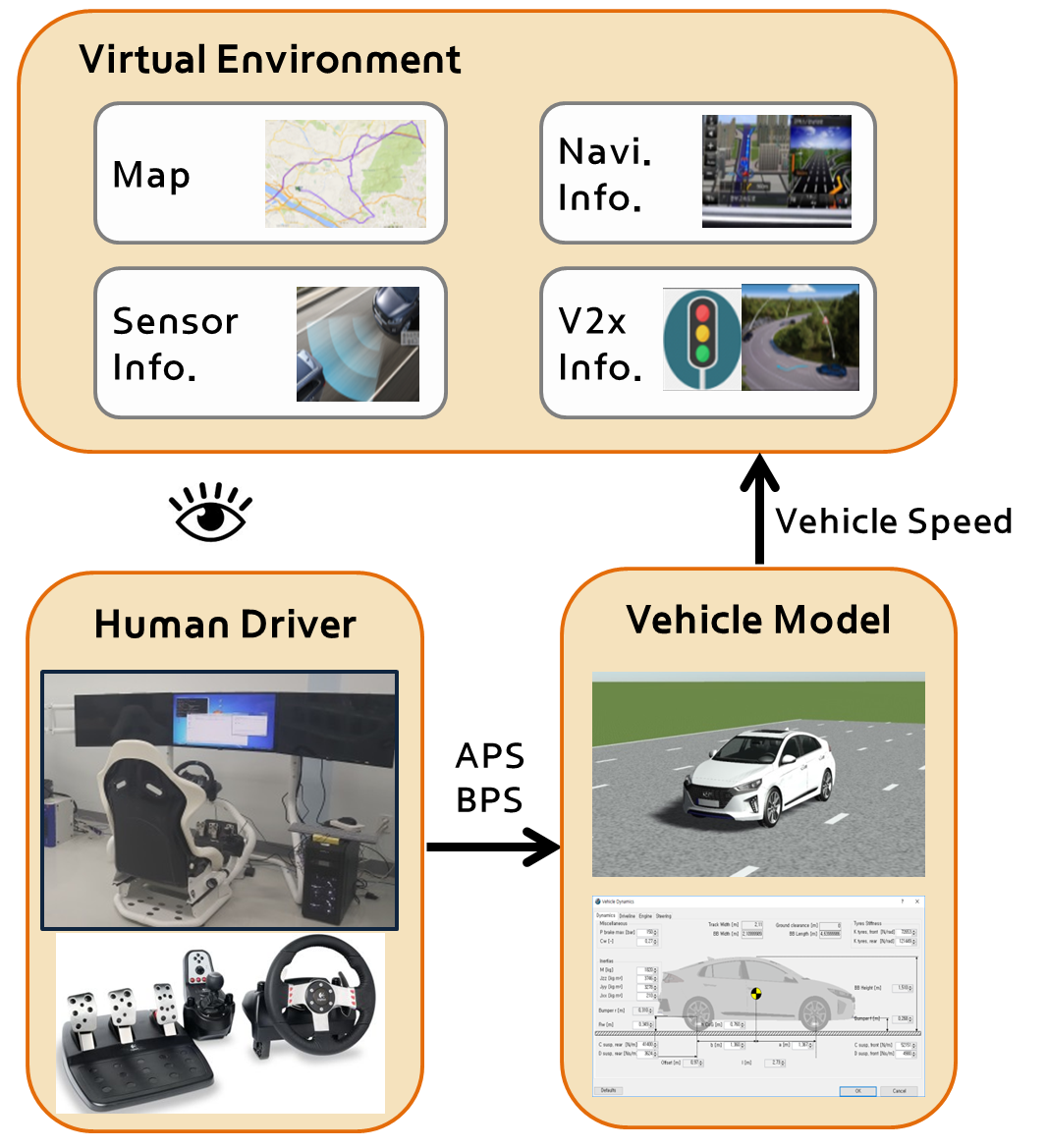}
	\caption{Illustrative configuration of the virtual driving test setup for the human-in-the-loop.}
	\label{fig:DILS}
\end{figure}
\begin{figure}[t]
	\centering
	\includegraphics[width=0.975\columnwidth]{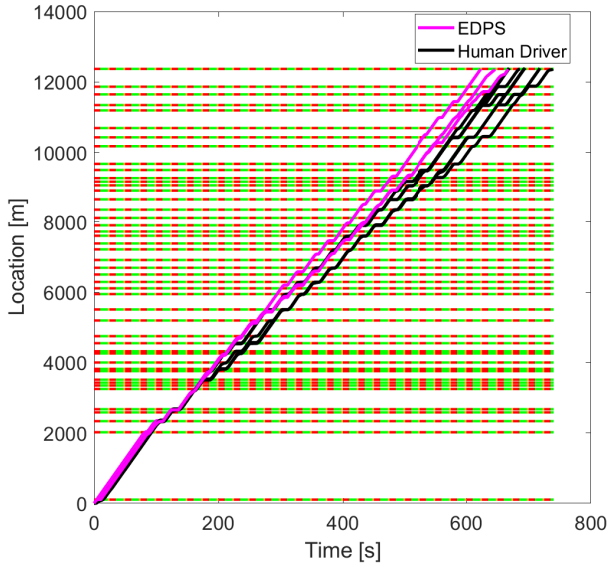}\vspace{-1mm}
	\caption{Comparisons in location profiles over total driving time and traffic SPaT: EDPS using three different preview distances (100, 150, 200 m) vs. Human drivers.}
	\label{fig:All_loc_time}
\end{figure}
\begin{figure}[t]
	\centering
	\includegraphics[width=0.975\columnwidth]{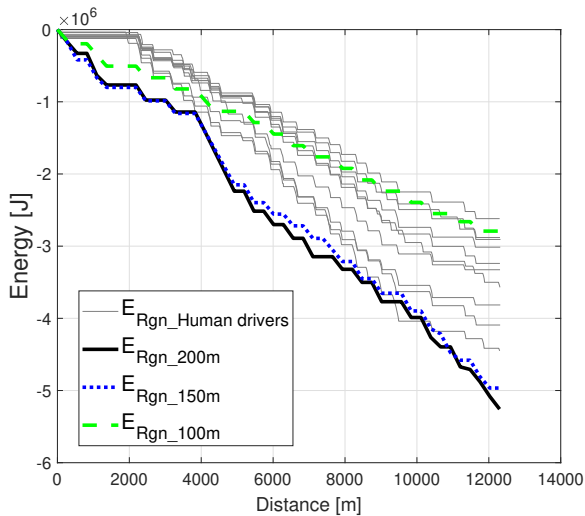}
	\caption{Comparisons in recuperated energy profiles over trip distances: EDPS using three different preview distances (100, 150, 200 m) vs. Human drivers.}
	\label{fig:Nrg_profs}
\end{figure}
\begin{figure}[t]
	\centering
	\includegraphics[width=0.975\columnwidth]{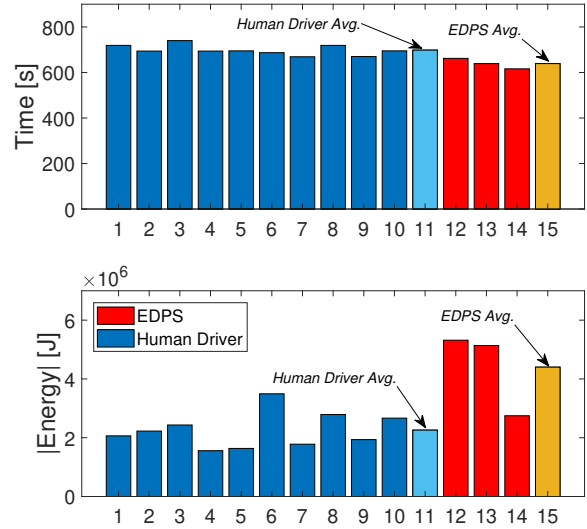}
	\caption{Comparisons in average recuperated energy and trip time: EDPS vs. Human drivers.}
	\label{fig:bar_DILS}
\end{figure}
\begin{figure}[t]
	\centering
	\includegraphics[width=0.975\columnwidth]{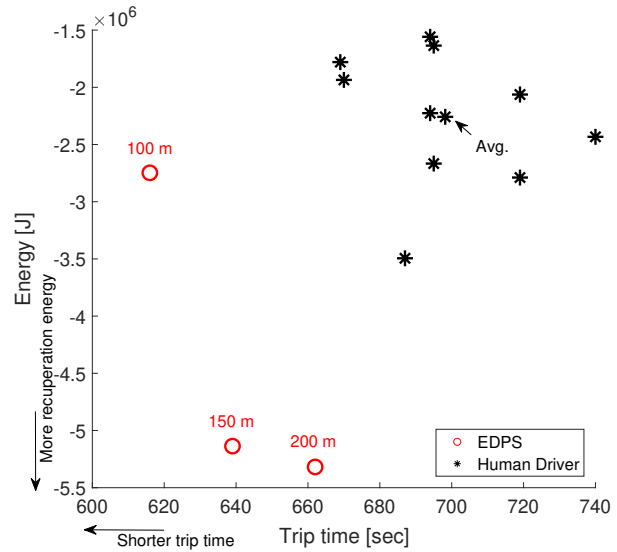}
	\caption{Comparisons in recuperated energy and trip time: EDPS using three different preview distances (100, 150, 200 m) vs. Human drivers.}
	\label{fig:pallette}
\end{figure}
%

\subsection{Comparison of EDPS to Human Drivers}
In principle, the concept of EDPS can be regarded as a semi-autonomous braking control that can be classified as levels $2$ or $3$. The acceleration is dependent on the decision of the driver while deceleration is automatically executed by energy-optimal deceleration control utilizing the look-ahead information when the deceleration event is detected beforehand.

To compare the deceleration performance of EDPS and a human driver in terms of energy recovery and basic deceleration performance, a driving test environment that simulates driving conditions for a human driver on the virtual route was established, as seen in Fig.~\ref{fig:DILS}. Each human driver operated the test driving equipment and recognized the traffic lights depending only on their vision, and the vehicle speed was limited to $90{\rm km/h}$ for a fair comparison with EDPS. The recuperated energy for human drivers is computed in a similar manner to the EDPS energy recovery calculation at each deceleration case.
The regenerative braking forces for human drivers are also constrained by considering the regeneration limit, as seen in~\eqref{F_rgn}.
The total braking forces are determined by considering brake pedal depths as seen in Fig.~\ref{fig:dec_bps_plot}. The data set was collected from the test-driving results obtained from $10$ human drivers in a virtual environment.

Comparing the test results of EDPS and human drivers, Fig.~\ref{fig:All_loc_time} and the top of Fig.~\ref{fig:bar_DILS} show that EDPS generally takes a shorter driving time compared with human drivers. The average results of $11^{\rm th}$ and $15^{\rm th}$ bar in Fig.~\ref{fig:bar_DILS} display that EDPS results in average time savings of $59\: {\rm s}$ in the total driving time. Since, human drivers normally consider the current traffic light status within their vision range, they decelerate unnecessarily by responding to the current red light, even though the traffic light will be changed to green light within a short time.

Regarding the accumulated energy of all human drivers (gray thin lines) and the EDPS accumulated energy with three preview distances, $200\:{\rm m}$, $150\:{\rm m}$, and $100\:{\rm m}$, Fig.~\ref{fig:Nrg_profs} indicates that the EDPS energy profiles with the preview distance of $200\:{\rm m}$ and $150\:{\rm m}$ can acquire more regenerative energy  compared with braking performed by a human driver because the preview information and the optimization process increase the regenerative braking force up to the regenerating limit of the electrified powertrain. The EDPS profile with preview distance of $100\:{\rm m}$ has a similar scale to that of the energy profile of human drivers.

The energy recovery results of Fig.~\ref{fig:pallette} and the bottom of Fig.~\ref{fig:bar_DILS} also exhibit greater benefits of using EDPS compared with those of human drivers because the recuperated energy of EDPS is twice as much as that of human drivers when comparing the $11^{\rm th}$ bar (the mean result of human drivers) with $15^{\rm th}$ bar (the mean result of EDPS). EDPS results in more recuperated energy than human drivers because it uses remote communication with exterior sources for taking a further preview distance and also to plan a more energy-efficient speed profile by the optimization scheme using the map information and traffic SPaT. In particular, the $200\:{\rm m}$ test case with the longest preview distance indicates the best energy-recuperation performance. In addition, while the $100\:{\rm m}$ test case records a shorter driving time than that of human drivers, the recuperated energy is the least among other EDPS test cases. This implies that the short preview distance decreases the opportunity for it to be regenerated owing to the shorter deceleration time compared with that of the other EDPS tests. Another interesting observation is that the energy recovery result for the $100\:{\rm m}$ case is similar to the averaged one for human drivers. For planning with consideration of energy consumption as well as energy recovery, $100\:{\rm m}$ preview might outperform human drivers in energy-saving, because it can provide opportunities of reducing the total duration of energy consumption.

\section{Conclusions}\label{sec5:conclusion}
This paper presents an optimal deceleration planning system called EDPS that provides a schedule of deceleration speed over a period of time when approaching an upcoming deceleration event. The objective is to maximize energy recuperation of regenerative braking and the proposed EDPS incorporates preview information available from CAV technology to determine optimal deceleration parameters. An optimal control problem is formulated with practical consideration of the energy recuperation model for the commercial $P_2$-type electrified vehicle. For energy-optimal regenerative braking strategy, a parameterized deceleration model is used and a feasible set of deceleration times of the model are sought to maximize the energy recuperation during the planning time. The deceleration time is the key design parameter to be determined and is constrained within an envelope reflecting the realistic deceleration limits obtained by dedicated real-world driving tests. In addition to the deceleration command limits, the system state constraints are dynamically adjusted by considering the   road load forces and the deceleration preference. For validation and verification, the proposed EDPS with different preview distances are tested and compared with human drivers in a facility of human-in-the-loop driving tests. In simulation case studies, the EDPS with longer preview distance attains a better energy-recuperated potential because the energy-optimal speed profile considering the system recuperation limit can be continued for a longer duration. In DILS tests, compared to 10 human drivers, the EDPS achieves a better energy-recuperation performance and a shorter trip time in virtual driving tests.

\section*{Acknowledgment}
The first and second authors are supported by Hyundai Motor Company and the third author (corresponding author) is partially supported by INHA UNIVERSITY Research Grant (INHA-60703).

\bibliographystyle{cas-model2-names}
\bibliography{edps}

\appendix
\subsection{Derivation of a Polynomial-based Deceleration Model}\label{sec:appendix}

This section presents the procedures for obtaining the main parameters of the deceleration model in~\eqref{dip}.
To find the parameters, consider a continuous deceleration profile model as
\begin{equation}\label{c_dip}
a_{p}\left(t\right)=r\alpha \theta\left(t\right)\left(1-\theta^{p}\left(t\right)\right)^{2},
\end{equation}
where
\begin{equation}\label{c_th}
\theta\left(t\right)=\frac{t}{t_{d}}.    
\end{equation}
To find the minimum of the deceleration profile model of (\ref{c_dip}), the derivative of (\ref{c_dip}) at $t=t_{M}$ defined as the time at minimum deceleration is written as
\begin{equation}\label{d_dip}
\begin{split}
\left. \frac{da_{d}}{dt} \right|_{t=t_{M}}
&=\frac{r\alpha}{t_{d}}\left(1-\left(2p+2\right)\theta_{M}^{p}+\left(2p+1\right)\theta_{M}^{2p}\right)\\
&=0,
\end{split}
\end{equation}
where $\theta_{M}=\theta\left(t_{M}\right)$.
Solving the quadratic equation in (\ref{d_dip}) yields the following solution set:
\begin{equation}\label{th_m}
\theta_{M}^{p}=\left\{ \begin{array}{ll}
\frac{1}{2p+1}, & \quad p \geq 0\\
\frac{-2p+1}{2p+1}, & \quad -0.5<p<0
\end{array}\right.
\end{equation}
In addition, using $a\left(t_{M}\right)=\alpha$, $r$ of (\ref{c_dip}) is expressed as
\begin{equation}\label{r}
r=\frac{1}{\theta_{M}\left(1-\theta_{M}^{p}\right)^{2}}.    
\end{equation}
With the substitution of (\ref{th_m}) for (\ref{r}), $r$ is presented as
\begin{equation*}
r=\left\{ \begin{array}{ll}
\frac{\left(1+2p\right)^{\frac{1}{p}+2}}{4p^{2}}, & \quad p \geq 0\\
\frac{\left(1+2p\right)^{\frac{1}{p}+2}}{16p^{2}\left(1-2p\right)^{\frac{1}{p}}}. & \quad -0.5<p<0
\end{array}\right.
\end{equation*}
Integrating (\ref{c_dip}),
\begin{equation}\label{v}
v_{f}=v_{i}+t_{d}r\alpha\theta^{2}\left(t\right)\left(0.5-2\frac{\theta^{p}\left(t\right)}{p+2}+\frac{\theta^{2p}\left(t\right)}{2p+2}\right).
\end{equation}
Using the definition of (\ref{c_th}) at $t=t_{d}$, $\theta\left(t_{d}\right)=1$ and (\ref{v}) are arranged as
\begin{equation*}
\frac{v_{f}-v_{i}}{t_{d}}=r\alpha q,
\end{equation*}
where $v_{f}=v\left(t_{d}\right)$ and $q$ is given by
\begin{equation*}
q=\frac{p^{2}}{\left(2p+2\right)\left(p+2\right)}.
\end{equation*}
Integrating (\ref{v}) from $t_{i}$ to $t_{d}$,
\begin{equation}\label{x}
x\left(t_{d}\right)-x\left(t_{i}\right)=v_{i}t_{d}+t_{d}^{2}r\alpha s,
\end{equation}
where
\begin{equation*}
s=\frac{1}{6}-\frac{2}{\left(p+2\right)\left(p+3\right)}+\frac{1}{\left(2p+2\right)\left(2p+3\right)}.
\end{equation*}
The ratio of the average speed during deceleration to the final speed, $\phi$, is defined as
\begin{equation}\label{phi}
\phi=\frac{v_{d}-v_{f}}{v_{i}-v_{f}},
\end{equation}
where
\begin{equation}\label{vd}
v_{d}=\frac{x\left(t_{d}\right)-x\left(t_{i}\right)}{t_{d}}.
\end{equation}
By substituting (\ref{x}) and (\ref{vd}) into (\ref{phi}), and performing some algebraic calculations, $\phi$ is rewritten as
\begin{equation}\label{phi2}
\phi=\frac{s}{q}=\frac{2p^{2}+15p+19}{3\left(p+3\right)\left(2p+3\right)}.
\end{equation}
Because the initial and final vehicle speeds during the planning time, $v_{i}$ and $v_{f}$, are given with the determination of an anticipative upcoming deceleration event, various candidates of $t_{d}$ ($N_{d}$ in the discrete domain) can make various values of $\phi$ from~\eqref{phi} and \eqref{vd}, which determines the values of $p$ from~\eqref{phi2}. With these candidate values of the parameter $p$, various deceleration profile candidates can be generated by the model~\eqref{c_dip}, as shown in Figs.~\ref{fig:Decel_cand} and \ref{fig:adp_decel}. 

\end{document}